\documentclass[11pt]{article}

\usepackage{caption}
\usepackage{booktabs} 
\usepackage[T1]{fontenc}
\usepackage[utf8]{inputenc}
\usepackage{lmodern}
\usepackage{float} 
\usepackage[margin=1in]{geometry}

\usepackage{amsmath,amssymb,amsthm,mathtools}
\usepackage{mathrsfs}
\usepackage{graphicx}
\usepackage{enumitem}
\usepackage{hyperref}
\usepackage[nameinlink,noabbrev]{cleveref}

\usepackage[all,cmtip]{xy}
\usepackage{pgfplots}

\hypersetup{
    colorlinks=true,
    linkcolor=black,
    citecolor=black,
    urlcolor=blue
}

\newtheorem{theorem}{Theorem}[section]
\newtheorem{lemma}[theorem]{Lemma}
\newtheorem{proposition}[theorem]{Proposition}
\newtheorem{corollary}[theorem]{Corollary}

\theoremstyle{definition}
\newtheorem{definition}[theorem]{Definition}
\newtheorem{remark}[theorem]{Remark}

\DeclareRobustCommand{\Dodot}{\ensuremath{\bodot{2}}}

\DeclareRobustCommand{\bodot}[1]{\ensuremath{\mathbin{ \odot_\text{\hspace{-0.5 pt}\scalebox{.775}{$[\! \:#1\! \:]$}}}}}

\title{\textbf{Irreversible $k$-Threshold Dynamics on\\
  Corona and Base-$b$ Corona Product Graphs}}
\author{Eric J. Moon, Soumya Bhoumik, Paul Flesher\\
Department of Mathematics\\
Fort Hays State University
}
\date{\today}

\begin{document}
\maketitle

\definecolor{light_gray}{gray}{0.7}

\begin{abstract}
We study the irreversible $k$-threshold process on corona-type graph products, including the corona product, the double corona product, and a generalized base-$b$ corona construction. Exact results are obtained for the irreversible $k$-threshold conversion number on corona and double corona product graphs through reduction lemmas that relate these graphs to smaller corona-type instances and to classical base graphs. The base-$b$ construction provides a unified framework extending the corona and double corona cases. A probabilistic refinement is also considered by studying the likelihood that a uniformly chosen minimum seed set yields complete activation. These results show how layered attachment structure influences both deterministic conversion behavior and probabilistic saturation in corona-type graph products.
\end{abstract}

\section{Introduction}

The irreversible $k$-threshold process is a discrete spreading model in which a vertex becomes active once at least $k$ of its neighbors are active, and remains active thereafter. A central problem is to determine the minimum size of an initial seed set that guarantees complete activation.

Since its introduction by Dreyer and Roberts \cite{dreyer2009}, the irreversible $k$-threshold process has been studied on a variety of classical graph families and graph products; see, for example, \cite{adams2015, shaheen2022, shaheen2022paths, abiad2024}. In contrast, corona-type constructions have received comparatively little attention in this setting, despite their layered structure and pronounced degree asymmetry.

Irreversible $k$-threshold conversion sets are studied on corona product graphs of the form $C_n \odot K_p$, together with related double corona and base-$b$ corona constructions. Exact threshold results are obtained through reduction arguments, and probabilistic behavior is also considered through random initial seed sets.

\section{Preliminaries}
\subsection{The Irreversible $k$-threshold Conversion Process}\label{Ik}
The irreversible $k$-threshold conversion process on a graph $G=(V,E)$ begins with an initial vertex set $S_0 \subseteq V(G)$. For each $t \ge 1$, the set $S_t$ is obtained from $S_{t-1}$ by adjoining every vertex with at least $k$ neighbors in $S_{t-1}$. Once a vertex becomes active, it remains active thereafter.

The set $S_0$ is called a \emph{seed set}. If $S_t = V(G)$ for some finite time $t$, then $S_0$ is an \emph{irreversible $k$-threshold conversion set} (I$k$CS) of $G$. The minimum cardinality of such a set is the \emph{irreversible $k$-threshold conversion number} of $G$, denoted $C_k(G)$.

A graph $G$ is \emph{$k$-inconvertible} if no proper subset of $V(G)$ is an I$k$CS; equivalently, if $C_k(G)=|V(G)|$.

\subsection{Basic Properties and Degree Constraints}
\begin{remark}\label{deg}
      Let \( G = (V, E) \) be a graph, and let \( S_0 \subset V \) be an I$k$CS, then the set \(U\) of vertices in \( G \) with degree less than \( k \) is contained in \( S_0 \); that is,
\[U=\{v \in V : \deg(v) < k\} \subseteq S_0\subset V.\]
\end{remark}

\begin{proposition}\label{k=1 spread}
For a graph $G$ where $G$ does not have isolated vertices, $C_1
(G)=1$    
\end{proposition}
\begin{proof}
Since $G$ has no isolated vertices, every vertex has at least one neighbor. If a single vertex is initially colored, then all of its neighbors become colored in the next step, and this process continues until the entire graph is colored. Hence $C_1(G)=1$.
\end{proof}

\begin{proposition}\label{prop-cycle}
    \[C_k(C_n) = 
    \begin{cases}
        1 & \text{if } k=1, \\
        \lceil \frac{n}{2}\rceil & \text{if } k=2, \\
        k\text{-inconvertible} & \text{if } k \ge 3, \\
    \end{cases}
    \]
\end{proposition}
\begin{proof}
The cases \(k=1\) and \(k=2\) are established in \cite{dreyer2009}.
For \(k \ge 3\), each vertex of \(C_n\) has degree 2, and thus no vertex
can ever acquire \(k\) colored neighbors. Hence no irreversible
\(k\)-threshold conversion set exists.
\end{proof}

\subsection{The Corona Product}

\begin{definition}
The \emph{corona product} of two graphs $G$ and $H$ (denoted $G \odot H$) is the graph obtained by taking one copy of $G$ and $|V(G)|$ copies of $H$, and joining the $i^{\rm th}$ vertex of $G$ to every vertex of the $i^{\rm th}$ copy of $H$.
\end{definition}

For the corona product $C_n \odot K_p$, the vertex set is
\[
V=\{v_i \mid 1\le i\le n\}\cup\{u_i^j \mid 1\le i\le n,\ 1\le j\le p\},
\]
where $v_i$ denotes the vertices of $C_n$ (referred to as \emph{inner vertices}), and $u_i^j$ represents the
$j^{\text{th}}$ vertex of the $i^{\text{th}}$ copy of $K_p$ attached to the vertex $v_i$.

\begin{remark}\label{C-size}
Note that the order of the graph 
$C_n \odot K_p$, is $$ |V(C_n \odot K_p)| = n \left(p+1\right),$$
whereas the size is 
$$|E(C_n \odot K_p)| = n+ n\binom{p}{2}+np = n\left(\frac{p(p+1)}{2}+1\right). $$
\end{remark}

\begin{figure}[H]
    \begin{minipage}{0.25\textwidth}
        \centering \[ \xygraph{
!{<0cm,0cm>;<0cm,0.75 cm>:<0.75 cm,0cm>::}
!{(0,0);a(0)**{}?(0.75)}*{\scalebox{1.5}{$\circ$}}="a1"
!{(0,0);a(120)**{}?(0.75)}*{\scalebox{1.5}{$\circ$}}="a2" 
!{(0,0);a(240)**{}?(0.75)}*{\scalebox{1.5}{$\circ$}}="a3" 
\\
!{(0,0);a(0)**{}?(1.75)}*{\scalebox{1.5}{$\circ$}}="b11"
\\
!{(0,0);a(125)**{}?(1.75)}*{\scalebox{1.5}{$\circ$}}="b21"
\\
!{(0,0);a(235)**{}?(1.75)}*{\scalebox{1.5}{$\circ$}}="b31"
\\
"a1"-"a2"   "a2"-"a3"   "a3"-"a1"
\\
"a1"-"b11"
"a2"-"b21"
"a3"-"b31"
}
\]
\vspace{-10 pt}\caption{\(C_3 \odot K_1\)}
\end{minipage}
\hfill
\begin{minipage}{0.25\textwidth}
 \[ \xygraph{
!{<0cm,0cm>;<0cm,0.75 cm>:<0.75 cm,0cm>::}
!{(0,0);a(0)**{}?(1)}*{\scalebox{1.25}{$\circ$}}="a1" 
        !{(0,0);a(0)**{}?(1.75)}
!{(0,0);a(60)**{}?(1)}*{\scalebox{1.5}{$\circ$}}="a2" 
        !{(0,0);a(60)**{}?(1.75)}
!{(0,0);a(120)**{}?(1)}*{\scalebox{1.5}{$\circ$}}="a3" 
        !{(0,0);a(120)**{}?(1.75)}
!{(0,0);a(180)**{}?(1)}*{\scalebox{1.5}{$\circ$}}="a4"
        !{(0,0);a(180)**{}?(1.75)}
!{(0,0);a(240)**{}?(1)}*{\scalebox{1.5}{$\circ$}}="a5"
        !{(0,0);a(240)**{}?(1.75)}
!{(0,0);a(300)**{}?(1)}*{\scalebox{1.5}{$\circ$}}="a6"
        !{(0,0);a(300)**{}?(1.75)}
\\
!{(0,0);a(-10)**{}?(1.75)}*{\scalebox{1.5}{$\circ$}}="b11"
!{(0,0);a(10)**{}?(1.75)}*{\scalebox{1.5}{$\circ$}}="b12"
\\
!{(0,0);a(50)**{}?(1.75)}*{\scalebox{1.5}{$\circ$}}="b21"
!{(0,0);a(70)**{}?(1.75)}*{\scalebox{1.5}{$\circ$}}="b22"
\\
!{(0,0);a(110)**{}?(1.75)}*{\scalebox{1.5}{$\circ$}}="b31"
!{(0,0);a(130)**{}?(1.75)}*{\scalebox{1.5}{$\circ$}}="b32"
\\
!{(0,0);a(170)**{}?(1.75)}*{\scalebox{1.5}{$\circ$}}="b41"
!{(0,0);a(190)**{}?(1.75)}*{\scalebox{1.5}{$\circ$}}="b42"
\\
!{(0,0);a(230)**{}?(1.75)}*{\scalebox{1.5}{$\circ$}}="b51"
!{(0,0);a(250)**{}?(1.75)}*{\scalebox{1.5}{$\circ$}}="b52"
\\
!{(0,0);a(290)**{}?(1.75)}*{\scalebox{1.5}{$\circ$}}="b61"
!{(0,0);a(310)**{}?(1.75)}*{\scalebox{1.5}{$\circ$}}="b62"
\\ 
"a1"-"a2"   "a2"-"a3"   "a3"-"a4"   "a4"-"a5"   "a5"-"a6"   "a6"-"a1"
\\
"b11"-"b12" 
\\
"b21"-"b22"
\\
"b31"-"b32"
\\
"b41"-"b42"
\\
"b51"-"b52"
\\
"b61"-"b62"
\\
"a1"-"b11"  "a1"-"b12"
"a2"-"b21"  "a2"-"b22"
"a3"-"b31"  "a3"-"b32"
"a4"-"b41"  "a4"-"b42"
"a5"-"b51"  "a5"-"b52"
"a6"-"b61"  "a6"-"b62"
}
\]      
\vspace{-10 pt}\caption{ \(C_6 \odot K_2\)}
    \end{minipage}
\hfill
\vspace{10 pt}\begin{minipage}{0.35\textwidth}
    \centering \[ \xygraph{
!{<0cm,0cm>;<0cm,1 cm>:<1 cm,0cm>::}
!{(0,0);a(0)**{}?(1)}*{\scalebox{1.5}{$\circ$}}="a1"
!{(0,0);a(40)**{}?(1)}*{\scalebox{1.5}{$\circ$}}="a2" 
!{(0,0);a(80)**{}?(1)}*{\scalebox{1.5}{$\circ$}}="a3" 
!{(0,0);a(120)**{}?(1)}*{\scalebox{1.5}{$\circ$}}="a4" 
!{(0,0);a(160)**{}?(1)}*{\scalebox{1.5}{$\circ$}}="a5"
!{(0,0);a(200)**{}?(1)}*{\scalebox{1.5}{$\circ$}}="a6"
!{(0,0);a(240)**{}?(1)}*{\scalebox{1.5}{$\circ$}}="a7"
!{(0,0);a(280)**{}?(1)}*{\scalebox{1.5}{$\circ$}}="a8"
!{(0,0);a(320)**{}?(1)}*{\scalebox{1.5}{$\circ$}}="a9"
\\ 
!{(0,0);a(-10)**{}?(1.7)}*{\scalebox{1.5}{$\circ$}}="b11"
!{(0,0);a(10)**{}?(1.7)}*{\scalebox{1.5}{$\circ$}}="b12"
!{(0,0);a(-8)**{}?(2.3)}*{\scalebox{1.5}{$\circ$}}="b13"
!{(0,0);a(8)**{}?(2.3)}*{\scalebox{1.5}{$\circ$}}="b14"
\\ 
!{(0,0);a(30)**{}?(1.7)}*{\scalebox{1.5}{$\circ$}}="b21"
!{(0,0);a(50)**{}?(1.7)}*{\scalebox{1.5}{$\circ$}}="b22"
!{(0,0);a(32)**{}?(2.3)}*{\scalebox{1.5}{$\circ$}}="b23"
!{(0,0);a(48)**{}?(2.3)}*{\scalebox{1.5}{$\circ$}}="b24"
\\ 
!{(0,0);a(70)**{}?(1.7)}*{\scalebox{1.5}{$\circ$}}="b31"
!{(0,0);a(90)**{}?(1.7)}*{\scalebox{1.5}{$\circ$}}="b32"
!{(0,0);a(72)**{}?(2.3)}*{\scalebox{1.5}{$\circ$}}="b33"
!{(0,0);a(88)**{}?(2.3)}*{\scalebox{1.5}{$\circ$}}="b34"
\\ 
!{(0,0);a(110)**{}?(1.7)}*{\scalebox{1.5}{$\circ$}}="b41"
!{(0,0);a(130)**{}?(1.7)}*{\scalebox{1.5}{$\circ$}}="b42"
!{(0,0);a(112)**{}?(2.3)}*{\scalebox{1.5}{$\circ$}}="b43"
!{(0,0);a(128)**{}?(2.3)}*{\scalebox{1.5}{$\circ$}}="b44"
\\ 
!{(0,0);a(150)**{}?(1.7)}*{\scalebox{1.5}{$\circ$}}="b51"
!{(0,0);a(170)**{}?(1.7)}*{\scalebox{1.5}{$\circ$}}="b52"
!{(0,0);a(152)**{}?(2.3)}*{\scalebox{1.5}{$\circ$}}="b53"
!{(0,0);a(168)**{}?(2.3)}*{\scalebox{1.5}{$\circ$}}="b54"
\\ 
!{(0,0);a(190)**{}?(1.7)}*{\scalebox{1.5}{$\circ$}}="b61"
!{(0,0);a(210)**{}?(1.7)}*{\scalebox{1.5}{$\circ$}}="b62"
!{(0,0);a(192)**{}?(2.3)}*{\scalebox{1.5}{$\circ$}}="b63"
!{(0,0);a(208)**{}?(2.3)}*{\scalebox{1.5}{$\circ$}}="b64"
\\ 
!{(0,0);a(230)**{}?(1.7)}*{\scalebox{1.5}{$\circ$}}="b71"
!{(0,0);a(250)**{}?(1.7)}*{\scalebox{1.5}{$\circ$}}="b72"
!{(0,0);a(232)**{}?(2.3)}*{\scalebox{1.5}{$\circ$}}="b73"
!{(0,0);a(248)**{}?(2.3)}*{\scalebox{1.5}{$\circ$}}="b74"
\\ 
!{(0,0);a(270)**{}?(1.7)}*{\scalebox{1.5}{$\circ$}}="b81"
!{(0,0);a(290)**{}?(1.7)}*{\scalebox{1.5}{$\circ$}}="b82"
!{(0,0);a(272)**{}?(2.3)}*{\scalebox{1.5}{$\circ$}}="b83"
!{(0,0);a(288)**{}?(2.3)}*{\scalebox{1.5}{$\circ$}}="b84"
\\ 
!{(0,0);a(310)**{}?(1.7)}*{\scalebox{1.5}{$\circ$}}="b91"
!{(0,0);a(330)**{}?(1.7)}*{\scalebox{1.5}{$\circ$}}="b92"
!{(0,0);a(312)**{}?(2.3)}*{\scalebox{1.5}{$\circ$}}="b93"
!{(0,0);a(328)**{}?(2.3)}*{\scalebox{1.5}{$\circ$}}="b94"
\\ 
"a1"-"a2"   "a2"-"a3"   "a3"-"a4"   "a4"-"a5"   "a5"-"a6"    "a6"-"a7"   "a7"-"a8"   "a8"-"a9"   "a9"-"a1"
\\ 
"b11"-"b12" "b11"-"b13" "b11"-"b14"
"b12"-"b13" "b12"-"b14"
"b13"-"b14" 
\\
"b21"-"b22" "b21"-"b23" "b21"-"b24" 
"b22"-"b23" "b22"-"b24" 
"b23"-"b24" 
\\
"b31"-"b32" "b31"-"b33" "b31"-"b34" 
"b32"-"b33" "b32"-"b34" 
"b33"-"b34"  
\\
"b41"-"b42" "b41"-"b43" "b41"-"b44"
"b42"-"b43" "b42"-"b44"
"b43"-"b44"
\\
"b51"-"b52" "b51"-"b53" "b51"-"b54"
"b52"-"b53" "b52"-"b54"
"b53"-"b54"  
\\
"b61"-"b62" "b61"-"b63" "b61"-"b64" 
"b62"-"b63" "b62"-"b64" 
"b63"-"b64" 
\\
"b71"-"b72" "b71"-"b73" "b71"-"b74"
"b72"-"b73" "b72"-"b74"
"b73"-"b74"
\\
"b81"-"b82" "b81"-"b83" "b81"-"b84"
"b82"-"b83" "b82"-"b84"
"b83"-"b84"
\\
"b91"-"b92" "b91"-"b93" "b91"-"b94"
"b92"-"b93" "b92"-"b94"
"b93"-"b94" 
\\ 
"a1"-"b11"  "a1"-"b12"  "a1"-"b13"  "a1"-"b14" 
"a2"-"b21"  "a2"-"b22"  "a2"-"b23"  "a2"-"b24" 
"a3"-"b31"  "a3"-"b32"  "a3"-"b33"  "a3"-"b34"  
"a4"-"b41"  "a4"-"b42"  "a4"-"b43"  "a4"-"b44" 
"a5"-"b51"  "a5"-"b52"  "a5"-"b53"  "a5"-"b54"  
"a6"-"b61"  "a6"-"b62"  "a6"-"b63"  "a6"-"b64"  
"a7"-"b71"  "a7"-"b72"  "a7"-"b73"  "a7"-"b74" 
"a8"-"b81"  "a8"-"b82"  "a8"-"b83"  "a8"-"b84"  
"a9"-"b91"  "a9"-"b92"  "a9"-"b93"  "a9"-"b94"
}
\]
\vspace{-8 pt}\caption{ \(C_9 \odot K_4\)}
  \end{minipage}
  \end{figure}
\subsection{The Double Corona Product}

\begin{definition}
The \emph{double corona product} of two graphs $G$ and $H$ (denoted $G \Dodot H$) is formed by taking two identical copies of $G$, denoted $G^{\prime}$ and $G^{\prime\prime}$, together with $g=|V(G)|$ copies of $H$. For each index $i \in \{1,2,\dots,g\}$, the $i^{\text{th}}$ vertex of $G^{\prime}$ and the $i^{\text{th}}$ vertex of $G^{\prime\prime}$ are joined to every vertex in the $i^{\text{th}}$ copy of $H$.

Suppose \(G\) has vertex set \( V(G) = \{v_1, v_2, \cdots, v_g\} \) and \( H\) has vertex set \( V(H) = \{u_1, u_2, \cdots, u_h\}\) then:
\begin{equation*}
    V(G \Dodot H) = V(G^{\prime}) \cup V(G^{\prime\prime}) \cup \bigcup_{i=1}^{g} V(H^{[i]})
\end{equation*}
\begin{equation*}
    E(G \Dodot H) = E(G^{\prime}) \cup E(G^{\prime\prime}) \cup \bigcup_{i=1}^{g} E(H^{[i]}) \cup \bigcup_{i=1}^{g} \left\{ (v_i^{\prime}, u),(v_i^{\prime\prime}, u) \mid \forall u \in V(H^{[i]}) \right\}
\end{equation*}

\end{definition}

For the double corona product $C_n\Dodot K_p$, the vertex set is
$$V=\{v_i\mid 1\le i\le n\} \cup \{u_i^j \mid 1\le i\le n, 1\le j\le p\} \cup \{w_i\mid 1\le i\le n\},$$
where $v_i$ denotes the vertices of $C_n^{\prime}$ (we refer to them as inner vertices), $w_i$ denotes the vertices of $C_n^{\prime\prime}$ (we refer to them as outer vertices), and where $u_i^j$ represents the $j^{\text{th}}$ vertex of the $i^{\text{th}}$ copy of the $K_p$ graph, which is connected with the $v_i$ vertex of $C_n^{\prime}$ and the $w_i$ vertex of $C_n^{\prime\prime}$.

\begin{remark}\label{D-size}
Note that the order of the graph 
$C_n \Dodot K_p$, is $$ |V(C_n \Dodot K_p)| = n \left(p+2\right),$$
whereas the size is 
$$|E(C_n \Dodot K_p)| = 2n+ n\binom{p}{2}+2np = n\left(\frac{(p+1)(p+2)}{2}+1\right). $$
\end{remark}

The $i^{\rm th}$ block of $C_n \Dodot K_p$ (denoted by $B^{[i]}$) is defined to be
\[
B^{[i]}=\{v_i,w_i\}\cup V(H^{[i]})
   =\{v_i,w_i\}\cup\{u_i^1,u_i^2,\dots,u_i^p\},
\]
and hence
\[
|B^{[i]}|=p+2.
\]

Blocks are classified according to the placement of seed vertices within them. 
A block is denoted
\begin{itemize}
\item $B$ (blank) if it contains no seed vertices;
\item $C$ (complete) if all vertices are seed or converted;
\item $I$ (inner) or $O$ (outer) if the unique seed vertex is an inner or outer vertex, respectively;
\item $M$ (middle) if exactly one vertex among $\{u_i^j \mid 1\le j\le p\}$ is a seed;
\item $T$ (two) if the block contains exactly two seed vertices, one inner and one outer.
\end{itemize}

\begin{figure}[H]
\begin{minipage}{0.475\textwidth}
\vspace{-53 pt}
\hspace{26 pt}
\rotatebox[]{45}{%
\hspace{-100 pt}
    \( \xygraph{
!{<0cm,0cm>;<0cm,1.25 cm>:<1.25 cm,0cm>::}
!{(0,0);a(0)**{}?(0.5)}*{\scalebox{1.25}{$\circ$}}="v1"
!{(0,0);a(90)**{}?(0.5)}*{\scalebox{1.25}{$\circ$}}="v2"
!{(0,0);a(180)**{}?(0.5)}*{\scalebox{1.25}{$\circ$}}="v3"
!{(0,0);a(270)**{}?(0.5)}*{\scalebox{1.25}{$\circ$}}="v4"
\\
!{(0,0);a(0)**{}?(2.5)}*{\scalebox{1.25}{$\circ$}}="w1"
!{(0,0);a(90)**{}?(2.5)}*{\scalebox{1.25}{$\circ$}}="w2"
!{(0,0);a(180)**{}?(2.5)}*{\scalebox{1.25}{$\circ$}}="w3"
!{(0,0);a(270)**{}?(2.5)}*{\scalebox{1.25}{$\circ$}}="w4"
\\ 
!{(0,0);a(0)**{}?(1.85)}*{\scalebox{1.25}{$\circ$}}="a1"
!{(0,0);a(14)**{}?(1.55)}*{\scalebox{1.25}{$\circ$}}="a2" 
!{(0,0);a(11)**{}?(1.1)}*{\scalebox{1.25}{$\circ$}}="a3" 
!{(0,0);a(-11)**{}?(1.1)}*{\scalebox{1.25}{$\circ$}}="a4" 
!{(0,0);a(-14)**{}?(1.55)}*{\scalebox{1.25}{$\circ$}}="a5"
\\ 
!{(0,0);a(90)**{}?(1.85)}*{\scalebox{1.25}{$\circ$}}="b1"
!{(0,0);a(104)**{}?(1.55)}*{\scalebox{1.25}{$\circ$}}="b2" 
!{(0,0);a(101)**{}?(1.1)}*{\scalebox{1.25}{$\circ$}}="b3" 
!{(0,0);a(79)**{}?(1.1)}*{\scalebox{1.25}{$\circ$}}="b4" 
!{(0,0);a(76)**{}?(1.55)}*{\scalebox{1.25}{$\circ$}}="b5"
\\ 
!{(0,0);a(180)**{}?(1.85)}*{\scalebox{1.25}{$\circ$}}="c1"
!{(0,0);a(194)**{}?(1.55)}*{\scalebox{1.25}{$\circ$}}="c2" 
!{(0,0);a(191)**{}?(1.1)}*{\scalebox{1.25}{$\circ$}}="c3" 
!{(0,0);a(169)**{}?(1.1)}*{\scalebox{1.25}{$\circ$}}="c4" 
!{(0,0);a(166)**{}?(1.55)}*{\scalebox{1.25}{$\circ$}}="c5"
\\ 
!{(0,0);a(270)**{}?(1.85)}*{\scalebox{1.25}{$\circ$}}="d1"
!{(0,0);a(284)**{}?(1.55)}*{\scalebox{1.25}{$\circ$}}="d2" 
!{(0,0);a(281)**{}?(1.1)}*{\scalebox{1.25}{$\circ$}}="d3" 
!{(0,0);a(259)**{}?(1.1)}*{\scalebox{1.25}{$\circ$}}="d4" 
!{(0,0);a(256)**{}?(1.55)}*{\scalebox{1.25}{$\circ$}}="d5"
\\
"v1"-"v2"   "v2"-"v3"   "v3"-"v4"   "v4"-"v1"
\\
"w1"-"w2"   "w2"-"w3"   "w3"-"w4"   "w4"-"w1"
\\
"a1"-"a2" "a1"-"a3" "a1"-"a4" "a1"-"a5"
"a2"-"a3" "a2"-"a4" "a2"-"a5"
"a3"-"a4" "a3"-"a5"
"a4"-"a5"
\\
"b1"-"b2" "b1"-"b3" "b1"-"b4" "b1"-"b5"
"b2"-"b3" "b2"-"b4" "b2"-"b5"
"b3"-"b4" "b3"-"b5"
"b4"-"b5"
\\
"c1"-"c2" "c1"-"c3" "c1"-"c4" "c1"-"c5"
"c2"-"c3" "c2"-"c4" "c2"-"c5"
"c3"-"c4" "c3"-"c5"
"c4"-"c5"
\\
"d1"-"d2" "d1"-"d3" "d1"-"d4" "d1"-"d5"
"d2"-"d3" "d2"-"d4" "d2"-"d5"
"d3"-"d4" "d3"-"d5"
"d4"-"d5"
\\
"v1"-"a1"  "v1"-@/^-0.17cm/"a2"  "v1"-@/^-0.0225cm/"a3" "v1"-@/_-0.0225cm/"a4"   "v1"-@/_-0.17cm/"a5"  
"v2"-"b1"  "v2"-@/^-0.17cm/"b2"  "v2"-@/^-0.0225cm/"b3" "v2"-@/_-0.0225cm/"b4"   "v2"-@/_-0.17cm/"b5"
"v3"-"c1"  "v3"-@/^-0.17cm/"c2"  "v3"-@/^-0.0225cm/"c3" "v3"-@/_-0.0225cm/"c4"   "v3"-@/_-0.17cm/"c5"
"v4"-"d1"  "v4"-@/^-0.17cm/"d2"  "v4"-@/^-0.0225cm/"d3" "v4"-@/_-0.0225cm/"d4"   "v4"-@/_-0.17cm/"d5"
\\
"w1"-"a1"  "w1"-@/_-0.12cm/"a2"  "w1"-@/_-0.09cm/"a3" "w1"-@/^-0.09cm/"a4"   "w1"-@/^-0.12cm/"a5"  
"w2"-"b1"  "w2"-@/_-0.12cm/"b2"  "w2"-@/_-0.09cm/"b3" "w2"-@/^-0.09cm/"b4"   "w2"-@/^-0.12cm/"b5"
"w3"-"c1"  "w3"-@/_-0.12cm/"c2"  "w3"-@/_-0.09cm/"c3" "w3"-@/^-0.09cm/"c4"   "w3"-@/^-0.12cm/"c5"
"w4"-"d1"  "w4"-@/_-0.12cm/"d2"  "w4"-@/_-0.09cm/"d3" "w4"-@/^-0.09cm/"d4"   "w4"-@/^-0.12cm/"d5"
    }
    \)
    }
    \vspace{3 pt}\caption{\(C_4 \Dodot K_5\)}
    \label{C4DK5}
\end{minipage}
\hfill
\begin{minipage}{0.475\textwidth}
    \centering \[ \xygraph{
!{<0cm,0cm>;<0cm,1 cm>:<-1 cm,0cm>::}
!{(0,0);a(0)**{}?(1)}*{\scalebox{1.25}{$\circ$}}="a1"
!{(0,0);a(72)**{}?(1)}*{\scalebox{1.25}{$\circ$}}="a2" 
!{(0,0);a(144)**{}?(1)}*{\scalebox{1.25}{$\circ$}}="a3" 
!{(0,0);a(216)**{}?(1)}*{\scalebox{1.25}{$\circ$}}="a4" 
!{(0,0);a(288)**{}?(1)}*{\scalebox{1.25}{$\circ$}}="a5"
\\
!{(0,0);a(0)**{}?(2.5)}*{\scalebox{1.25}{$\circ$}}="c1"
!{(0,0);a(72)**{}?(2.5)}*{\scalebox{1.25}{$\circ$}}="c2" 
!{(0,0);a(144)**{}?(2.5)}*{\scalebox{1.25}{$\circ$}}="c3" 
!{(0,0);a(216)**{}?(2.5)}*{\scalebox{1.25}{$\circ$}}="c4" 
!{(0,0);a(288)**{}?(2.5)}*{\scalebox{1.25}{$\circ$}}="c5"
\\ 
!{(0,0);a(-10)**{}?(1.5)}*{\scalebox{1.25}{$\circ$}}="b11"
!{(0,0);a(10)**{}?(1.5)}*{\scalebox{1.25}{$\circ$}}="b12"
!{(0,0);a(0)**{}?(1.9)}*{\scalebox{1.25}{$\circ$}}="b13"
\\ 
!{(0,0);a(62)**{}?(1.5)}*{\scalebox{1.25}{$\circ$}}="b21"
!{(0,0);a(82)**{}?(1.5)}*{\scalebox{1.25}{$\circ$}}="b22"
!{(0,0);a(72)**{}?(1.9)}*{\scalebox{1.25}{$\circ$}}="b23"
\\ 
!{(0,0);a(134)**{}?(1.5)}*{\scalebox{1.25}{$\circ$}}="b31"
!{(0,0);a(154)**{}?(1.5)}*{\scalebox{1.25}{$\circ$}}="b32"
!{(0,0);a(144)**{}?(1.9)}*{\scalebox{1.25}{$\circ$}}="b33"
\\ 
!{(0,0);a(206)**{}?(1.5)}*{\scalebox{1.25}{$\circ$}}="b41"
!{(0,0);a(226)**{}?(1.5)}*{\scalebox{1.25}{$\circ$}}="b42"
!{(0,0);a(216)**{}?(1.9)}*{\scalebox{1.25}{$\circ$}}="b43"
\\ 
!{(0,0);a(278)**{}?(1.5)}*{\scalebox{1.25}{$\circ$}}="b51"
!{(0,0);a(298)**{}?(1.5)}*{\scalebox{1.25}{$\circ$}}="b52"
!{(0,0);a(288)**{}?(1.9)}*{\scalebox{1.25}{$\circ$}}="b53"
\\ 
"a1"-"a2"   "a2"-"a3"   "a3"-"a4"   "a4"-"a5"   "a5"-"a1"
\\ 
"c1"-"c2"   "c2"-"c3"   "c3"-"c4"   "c4"-"c5"   "c5"-"c1"
\\ 
"b11"-"b12" "b12"-"b13" "b13"-"b11" 
"b21"-"b22" "b22"-"b23" "b23"-"b21" 
"b31"-"b32" "b32"-"b33" "b33"-"b31" 
"b41"-"b42" "b42"-"b43" "b43"-"b41" 
"b51"-"b52" "b52"-"b53" "b53"-"b51" 
\\ 
"a1"-"b11"  "a1"-"b12"  "a1"-"b13"
"a2"-"b21"  "a2"-"b22"  "a2"-"b23"
"a3"-"b31"  "a3"-"b32"  "a3"-"b33"
"a4"-"b41"  "a4"-"b42"  "a4"-"b43"
"a5"-"b51"  "a5"-"b52"  "a4"-"b43"
\\ 
"c1"-"b11"  "c1"-"b12"  "c1"-"b13" 
"c2"-"b21"  "c2"-"b22"  "c2"-"b23" 
"c3"-"b31"  "c3"-"b32"  "c3"-"b33" 
"c4"-"b41"  "c4"-"b42"  "c4"-"b43" 
"c5"-"b51"  "c5"-"b52"  "c5"-"b53" 
}
\]
\vspace{-10 pt}\caption{\(C_5 \Dodot K_3\)}
\label{C5DK3}
\end{minipage}
\end{figure}

\begin{remark}\label{order/size}
From the definitions of corona and double corona graphs, it is clear that the graphs \(C_n \odot K_p\) and \(C_n  \Dodot  K_{p-1}\) have the same number of vertices and edges, but arise from different constructions. We are interested in how this difference will affect the irreversible \(k\)-threshold process, which is discussed in Section \ref{comparison}.
\end{remark}

\subsection{The Base-$b$ Corona Product}

The base-$b$ corona product generalizes the corona and double corona constructions to an arbitrary number of layers; the cases $b=1$ and $b=2$ recover the classical corona product and the double corona product, respectively.

\begin{definition}
The \emph{base-$b$ corona product} of two graphs $G$ and $H$ (denoted $G \bodot{b} H$) is formed by taking $b$ identical copies (or layers) of $G$, denoted $G^{[1]}, G^{[2]}, \dots, G^{[b]}$, together with $g=|V(G)|$ copies of $H$,
denoted $H^{[1]}, H^{[2]}, \dots, H^{[g]}$. For each layer $\ell \in \{1,2,\dots,b\}$ and each index $i \in \{1,2,\dots,g\}$, the $i^{\text{th}}$ vertex of $G^{[\ell]}$ is joined to every vertex in $H^{[i]}$ (the $i^{\text{th}}$ copy of~$H$).

Suppose \( G^{[\ell]} \) has vertex set \( V(G^{[\ell]}) = \{v_1^{[\ell]}, v_2^{[\ell]}, \cdots, v_g^{[\ell]}\} \) and \( H^{[i]} \) has vertex set \\ \( V(H^{[i]})=\{u_1^{[i]}, u_2^{[i]}, \cdots, u_h^{[i]}\}\) then:
\begin{equation*}
 V(G \bodot{b} H)  = \bigcup_{\ell=1}^{b} V(G^{[\ell]}) \cup \bigcup_{i=1}^{g} V(H^{[i]}) 
\end{equation*}
\begin{equation*}
 E(G \bodot{b} H) = \bigcup_{\ell=1}^{b} E(G^{[\ell]}) \cup \bigcup_{i=1}^{g} E(H^{[i]}) \cup \bigcup_{\ell=1}^{b}\bigcup_{i=1}^{g}\left\{(v_i^{[\ell]}, u) \mid \forall u \in V(H^{[i]}) \right\} 
\end{equation*}
\end{definition}
\begin{remark}
Note that the order of the graph 
$G \bodot{b} H$, is $$ |V(G \bodot{b} H)| = g \left(b+h\right),$$
whereas the size is 
$$|E(G \bodot{b} H)| = b\cdot \vert \!\:E(G) \!\:\vert + g \cdot \vert \!\:E(H) \!\:\vert + bgh. $$
 \end{remark}
The $i^{\rm th}$ zone of $G \bodot{b} H$ (denoted by $Z^{[i]}\:\!\!$) is defined to be:
    \[Z^{[i]} = \bigcup_{\ell=1}^{b} v_i^{[\ell]} \cup V(H^{[i]}) = \left\{v_i^{[1]}, v_i^{[2]}, \cdots, v_i^{[b]}\right\} \cup \left\{u_1^{[i]}, u_2^{[i]}, \cdots , u_h^{[i]} \right\},\]
    and hence
    \[\!\: \vert Z^{[i]} \!\: \vert = b+h\]
Zones play the same structural role in the base-$b$ corona product as blocks do in the double corona product.

\begin{remark}

The following special cases are immediate.
\begin{itemize}
  \item For $G = K_0$ (the empty graph), 
  \[K_0 \bodot{b} H \cong K_0.\]
  \item For $H = K_0$ (the empty graph), 
  \[
  G \bodot{b} K_0 \;\cong\; \bigcup_{\ell=1}^b G^{[\ell]},
  \]
  i.e., the union of $b$ disjoint copies of $G$.
  \item For $b=1$, we recover the classical corona product:
  \[
  G \bodot{1} H \;\cong\; G \odot H.
  \]
  \item For $b=2$, we obtain a distinct member of the family.
\end{itemize}
\label{special_bodot}
\end{remark}

In the case of $b=2$ (the double corona) $C_n \Dodot K_p$, we let $\{v_i^{[1]} \, | 1\le i\le n\} := \{v_i \, | 1\le i\le n\}$ and $\{v_i^{[2]} \, | 1\le i\le n\} := \{w_i \, | 1\le i\le n\}$.
\[
V(C_n \Dodot K_p)=\{v_i \mid 1\le i\le n\}\,\cup\,\{w_i \mid 1\le i\le n\}\,\cup\,\{u_j^{[i]} \mid 1\le i\le n,\; 1\le j\le p\},
\]
We refer to $I=\{v_1,\dots,v_n\}$ as the inner vertices and $O=\{w_1,\dots,w_n\}$ as the outer vertices.

\begin{proposition}\label{b-corona info}
Let $G=C_n$ and $H=K_p$. Then:
\begin{itemize}
\item \(\vert \!\:V(C_n \bodot{b} K_p) \!\: \vert = n \cdot\left(b+p\right) \)
\item \(\vert \!\:E(C_n \bodot{b} K_p) \!\: \vert = n\left(b+\binom{p}{2}+bp\right) = n\left(\frac{p(p+2b-1)}{2}+b\right) \)
\item \( \deg(v_i^{[\ell]}) = 2+p \)

\item If \(p\ge 1, \; \deg(u_j^{[i]}) = b+p-1 \)
\end{itemize}
\end{proposition}

\begin{remark}
Specializing Proposition~\ref{b-corona info} to $b=1$ and $b=2$ recovers the order and size formulas for $C_n \odot K_p$ and $C_n \Dodot K_p$ given earlier in Remarks \ref{C-size} and \ref{D-size} respectively.
\end{remark}

\begin{figure}[H]
    \centering
    \begin{minipage}{0.465\textwidth}
\includegraphics[height= 5.2 cm]{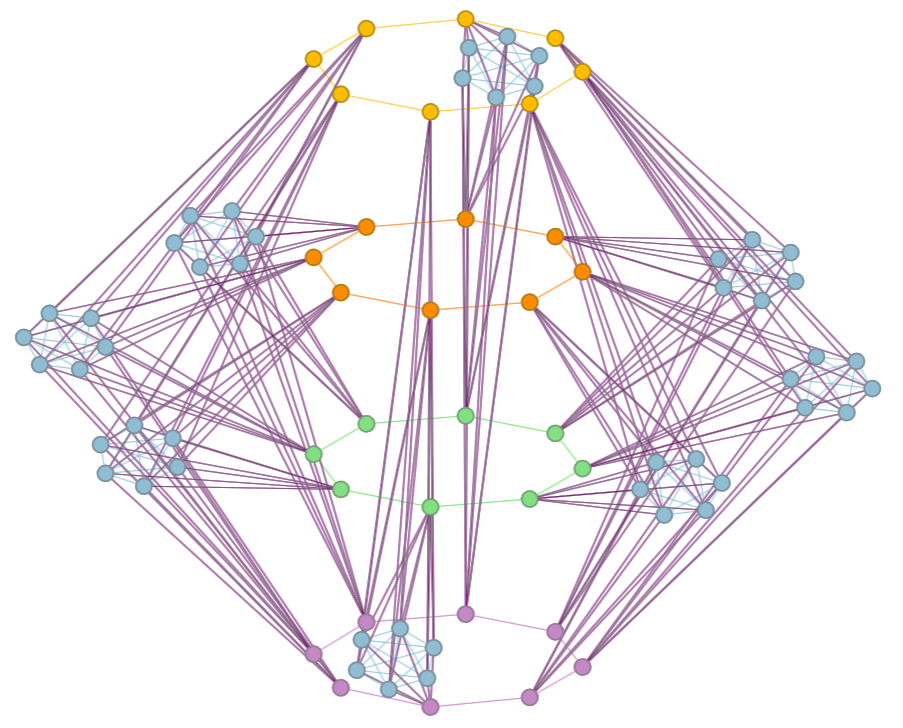}
        \caption{$C_8 \bodot{4} K_6$}
        \label{C8b4K6}
    \end{minipage}
    \hfill
    \begin{minipage}{0.515\textwidth}
\hspace{5 pt}\includegraphics[height= 7 cm]{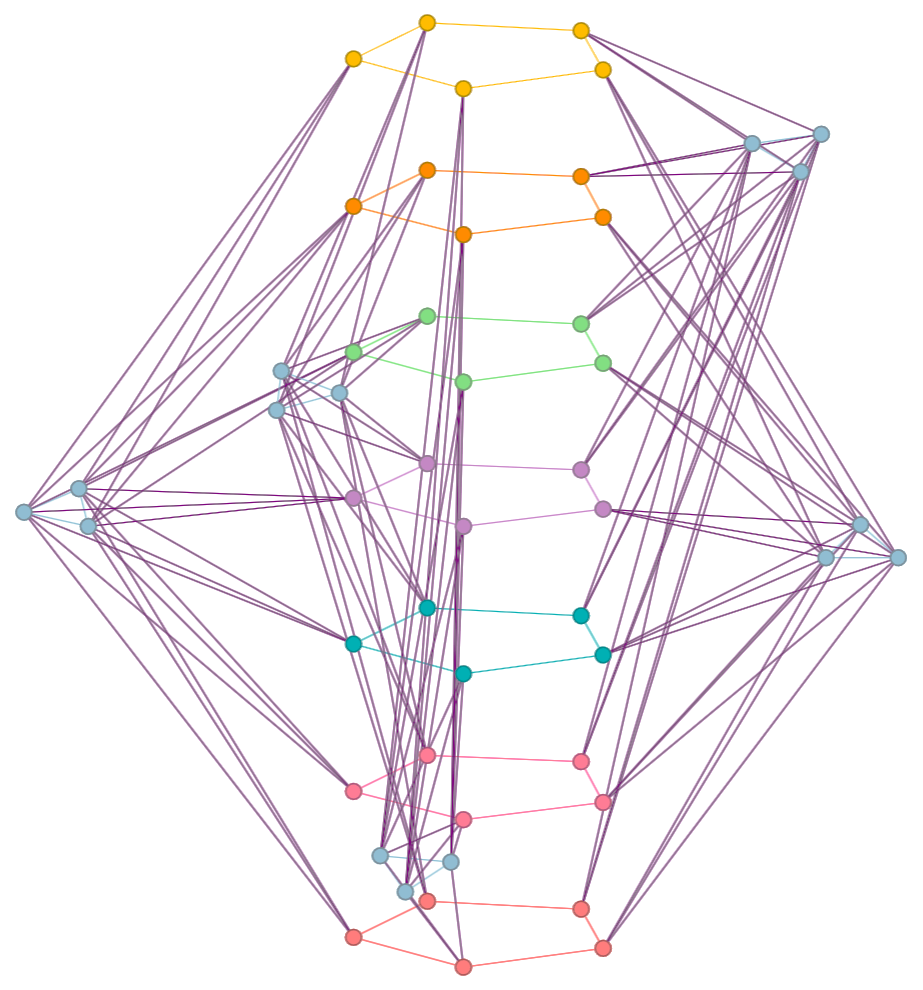}
        \caption{$C_5 \bodot{7} K_3$}

        \label{C5b7K3}
    \end{minipage}
    \end{figure}

\section{Threshold Dynamics on Corona Product Graphs}\label{C Chap}
The irreversible $k$-threshold process on corona product graphs admits a reduction to smaller instances, leading to exact threshold results for $C_n \odot K_p$.

\begin{lemma}\label{C reduct}
$C_k(C_n \odot K_p) = r\cdot n + C_{k-r}(C_n \odot K_{p-r})$ where $r=\min\{(k-1), p\}$.
\end{lemma}
\begin{proof}
If \(k \le p\), each complete graph in \(S_0\) must contain at least \(k-1\) seed vertices. If \( k > p \), all \( p \) vertices in each complete graph are necessary. Hence every vertex is adjacent to at least \(k-1\) colored vertices.

For \(k \le p\), remove exactly \(k-1\) colored vertices from each complete graph. The resulting graph $C_n\odot K_{p-(k-1)}$ then requires only a threshold of $1$ for completion.

If \(k>p\), all \(p\) vertices in each complete graph are removed, and the problem reduces to \(C_n \odot K_0\) with threshold \(k-p\). Thus
\[C_k(C_n \odot K_p) =
\begin{cases}
(k-1)n + C_{1}(C_n \odot K_{p - (k-1)}), & \text{if } k < p+1, \\
pn + C_{k - p}(C_n \odot K_0), & \text{if } k \ge p+1.
\end{cases} 
\]
\end{proof}

\begin{figure}[H]
    \centering
    \centering
  \begin{minipage}{0.245\textwidth}
    \centering
    \[
      \xygraph{
        !{<0cm,0cm>;<0cm,0.9cm>:<0.9cm,0cm>::}
        !{(0,0);a(0)**{}?(1)}*{\scalebox{1.25}{$\circ$}}="a1"
        !{(0,0);a(45)**{}?(1)}*{\scalebox{1.25}{$\circ$}}="a2"
        !{(0,0);a(90)**{}?(1)}*{\scalebox{1.25}{$\circ$}}="a3"
        !{(0,0);a(135)**{}?(1)}*{\scalebox{1.25}{$\circ$}}="a4"
        !{(0,0);a(180)**{}?(1)}*{\scalebox{1.25}{$\circ$}}="a5"
        !{(0,0);a(225)**{}?(1)}*{\scalebox{1.25}{$\circ$}}="a6"
        !{(0,0);a(270)**{}?(1)}*{\scalebox{1.25}{$\circ$}}="a7"
        !{(0,0);a(315)**{}?(1)}*{\scalebox{1.25}{$\circ$}}="a8"
        %
        !{(0,0);a(-10)**{}?(1.5)}*{\scalebox{1.25}{$\bullet$}}="b11"
        !{(0,0);a(10)**{}?(1.5)}*{\scalebox{1.25}{$\bullet$}}="b12"
        !{(0,0);a(0)**{}?(1.9)}*{\scalebox{1.25}{$\circ$}}="b13"
        !{(0,0);a(35)**{}?(1.5)}*{\scalebox{1.25}{$\bullet$}}="b21"
        !{(0,0);a(55)**{}?(1.5)}*{\scalebox{1.25}{$\bullet$}}="b22"
        !{(0,0);a(45)**{}?(1.9)}*{\scalebox{1.25}{$\circ$}}="b23"
        !{(0,0);a(80)**{}?(1.5)}*{\scalebox{1.25}{$\bullet$}}="b31"
        !{(0,0);a(100)**{}?(1.5)}*{\scalebox{1.25}{$\bullet$}}="b32"
        !{(0,0);a(90)**{}?(1.9)}*{\scalebox{1.25}{$\circ$}}="b33"
        !{(0,0);a(125)**{}?(1.5)}*{\scalebox{1.25}{$\bullet$}}="b41"
        !{(0,0);a(145)**{}?(1.5)}*{\scalebox{1.25}{$\bullet$}}="b42"
        !{(0,0);a(135)**{}?(1.9)}*{\scalebox{1.25}{$\circ$}}="b43"
        !{(0,0);a(170)**{}?(1.5)}*{\scalebox{1.25}{$\bullet$}}="b51"
        !{(0,0);a(190)**{}?(1.5)}*{\scalebox{1.25}{$\bullet$}}="b52"
        !{(0,0);a(180)**{}?(1.9)}*{\scalebox{1.25}{$\circ$}}="b53"
        !{(0,0);a(215)**{}?(1.5)}*{\scalebox{1.25}{$\bullet$}}="b61"
        !{(0,0);a(235)**{}?(1.5)}*{\scalebox{1.25}{$\bullet$}}="b62"
        !{(0,0);a(225)**{}?(1.9)}*{\scalebox{1.25}{$\circ$}}="b63"
        !{(0,0);a(260)**{}?(1.5)}*{\scalebox{1.25}{$\bullet$}}="b71"
        !{(0,0);a(280)**{}?(1.5)}*{\scalebox{1.25}{$\bullet$}}="b72"
        !{(0,0);a(270)**{}?(1.9)}*{\scalebox{1.25}{$\circ$}}="b73"
        !{(0,0);a(305)**{}?(1.5)}*{\scalebox{1.25}{$\bullet$}}="b81"
        !{(0,0);a(325)**{}?(1.5)}*{\scalebox{1.25}{$\bullet$}}="b82"
        !{(0,0);a(315)**{}?(1.9)}*{\scalebox{1.25}{$\circ$}}="b83"
        %
        "a1"-"a2" "a2"-"a3" "a3"-"a4" "a4"-"a5"
        "a5"-"a6" "a6"-"a7" "a7"-"a8" "a8"-"a1"
        "b11"-"b12" "b12"-"b13" "b13"-"b11"
        "b21"-"b22" "b22"-"b23" "b23"-"b21"
        "b31"-"b32" "b32"-"b33" "b33"-"b31"
        "b41"-"b42" "b42"-"b43" "b43"-"b41"
        "b51"-"b52" "b52"-"b53" "b53"-"b51"
        "b61"-"b62" "b62"-"b63" "b63"-"b61"
        "b71"-"b72" "b72"-"b73" "b73"-"b71"
        "b81"-"b82" "b82"-"b83" "b83"-"b81"
        "a1"-"b11" "a1"-"b12" "a1"-"b13"
        "a2"-"b21" "a2"-"b22" "a2"-"b23"
        "a3"-"b31" "a3"-"b32" "a3"-"b33"
        "a4"-"b41" "a4"-"b42" "a4"-"b43"
        "a5"-"b51" "a5"-"b52" "a5"-"b53"
        "a6"-"b61" "a6"-"b62" "a6"-"b63"
        "a7"-"b71" "a7"-"b72" "a7"-"b73"
        "a8"-"b81" "a8"-"b82" "a8"-"b83"
      }
    \]
    \vspace{-10 pt}\caption*{\(C_3(C_8\odot K_3)\)}
  \end{minipage}
  \hfill
  \begin{minipage}{0.245\textwidth}
    \centering
    \[
      \xygraph{
        !{<0cm,0cm>;<0cm,0.9cm>:<0.9cm,0cm>::}
        !{(0,0);a(0)**{}?(1)}*{\scalebox{1.25}{$\circ$}}="a1"
        !{(0,0);a(45)**{}?(1)}*{\scalebox{1.25}{$\circ$}}="a2"
        !{(0,0);a(90)**{}?(1)}*{\scalebox{1.25}{$\circ$}}="a3"
        !{(0,0);a(135)**{}?(1)}*{\scalebox{1.25}{$\circ$}}="a4"
        !{(0,0);a(180)**{}?(1)}*{\scalebox{1.25}{$\circ$}}="a5"
        !{(0,0);a(225)**{}?(1)}*{\scalebox{1.25}{$\circ$}}="a6"
        !{(0,0);a(270)**{}?(1)}*{\scalebox{1.25}{$\circ$}}="a7"
        !{(0,0);a(315)**{}?(1)}*{\scalebox{1.25}{$\circ$}}="a8"
        %
        !{(0,0);a(-10)**{}?(1.5)}*{\textcolor{light_gray}{\scalebox{1.05}{$\bullet$}}}="g11"
        !{(0,0);a(10)**{}?(1.5)}*{\textcolor{light_gray}{\scalebox{1.05}{$\bullet$}}}="g12"
        !{(0,0);a(0)**{}?(1.9)}*{\scalebox{1.25}{$\circ$}}="g13"
        !{(0,0);a(35)**{}?(1.5)}*{\textcolor{light_gray}{\scalebox{1.05}{$\bullet$}}}="g21"
        !{(0,0);a(55)**{}?(1.5)}*{\textcolor{light_gray}{\scalebox{1.05}{$\bullet$}}}="g22"
        !{(0,0);a(45)**{}?(1.9)}*{\scalebox{1.25}{$\circ$}}="g23"
        !{(0,0);a(80)**{}?(1.5)}*{\textcolor{light_gray}{\scalebox{1.05}{$\bullet$}}}="g31"
        !{(0,0);a(100)**{}?(1.5)}*{\textcolor{light_gray}{\scalebox{1.05}{$\bullet$}}}="g32"
        !{(0,0);a(90)**{}?(1.9)}*{\scalebox{1.25}{$\circ$}}="g33"
        !{(0,0);a(125)**{}?(1.5)}*{\textcolor{light_gray}{\scalebox{1.05}{$\bullet$}}}="g41"
        !{(0,0);a(145)**{}?(1.5)}*{\textcolor{light_gray}{\scalebox{1.05}{$\bullet$}}}="g42"
        !{(0,0);a(135)**{}?(1.9)}*{\scalebox{1.25}{$\circ$}}="g43"
        !{(0,0);a(170)**{}?(1.5)}*{\textcolor{light_gray}{\scalebox{1.05}{$\bullet$}}}="g51"
        !{(0,0);a(190)**{}?(1.5)}*{\textcolor{light_gray}{\scalebox{1.05}{$\bullet$}}}="g52"
        !{(0,0);a(180)**{}?(1.9)}*{\scalebox{1.25}{$\circ$}}="g53"
        !{(0,0);a(215)**{}?(1.5)}*{\textcolor{light_gray}{\scalebox{1.05}{$\bullet$}}}="g61"
        !{(0,0);a(235)**{}?(1.5)}*{\textcolor{light_gray}{\scalebox{1.05}{$\bullet$}}}="g62"
        !{(0,0);a(225)**{}?(1.9)}*{\scalebox{1.25}{$\circ$}}="g63"
        !{(0,0);a(260)**{}?(1.5)}*{\textcolor{light_gray}{\scalebox{1.05}{$\bullet$}}}="g71"
        !{(0,0);a(280)**{}?(1.5)}*{\textcolor{light_gray}{\scalebox{1.05}{$\bullet$}}}="g72"
        !{(0,0);a(270)**{}?(1.9)}*{\scalebox{1.25}{$\circ$}}="g73"
        !{(0,0);a(305)**{}?(1.5)}*{\textcolor{light_gray}{\scalebox{1.05}{$\bullet$}}}="g81"
        !{(0,0);a(325)**{}?(1.5)}*{\textcolor{light_gray}{\scalebox{1.05}{$\bullet$}}}="g82"
        !{(0,0);a(315)**{}?(1.9)}*{\scalebox{1.25}{$\circ$}}="g83"
        %
        "a1"-"a2" "a2"-"a3" "a3"-"a4" "a4"-"a5"
        "a5"-"a6" "a6"-"a7" "a7"-"a8" "a8"-"a1"
        "g11"-@{.}"g12" "g12"-@{.}"g13" "g13"-@{.}"g11"
        "g21"-@{.}"g22" "g22"-@{.}"g23" "g23"-@{.}"g21"
        "g31"-@{.}"g32" "g32"-@{.}"g33" "g33"-@{.}"g31"
        "g41"-@{.}"g42" "g42"-@{.}"g43" "g43"-@{.}"g41"
        "g51"-@{.}"g52" "g52"-@{.}"g53" "g53"-@{.}"g51"
        "g61"-@{.}"g62" "g62"-@{.}"g63" "g63"-@{.}"g61"
        "g71"-@{.}"g72" "g72"-@{.}"g73" "g73"-@{.}"g71"
        "g81"-@{.}"g82" "g82"-@{.}"g83" "g83"-@{.}"g81"
        "a1"-@{.}"g11" "a1"-@{.}"g12" "a1"-"g13"
        "a2"-@{.}"g21" "a2"-@{.}"g22" "a2"-"g23"
        "a3"-@{.}"g31" "a3"-@{.}"g32" "a3"-"g33"
        "a4"-@{.}"g41" "a4"-@{.}"g42" "a4"-"g43"
        "a5"-@{.}"g51" "a5"-@{.}"g52" "a5"-"g53"
        "a6"-@{.}"g61" "a6"-@{.}"g62" "a6"-"g63"
        "a7"-@{.}"g71" "a7"-@{.}"g72" "a7"-"g73"
        "a8"-@{.}"g81" "a8"-@{.}"g82" "a8"-"g83"
      }
    \]
    \vspace{-10 pt}\caption*{\(C_{3-2}(C_8 \odot K_{3-2})\)}
  \end{minipage}
  \hfill
  \begin{minipage}{0.245\textwidth}
    \centering
    \[
      \xygraph{
        !{<0cm,0cm>;<0cm,0.9cm>:<0.9cm,0cm>::}
        !{(0,0);a(0)**{}?(1)}*{\scalebox{1.5}{$\circ$}}="a1"
        !{(0,0);a(45)**{}?(1)}*{\scalebox{1.5}{$\circ$}}="a2"
        !{(0,0);a(90)**{}?(1)}*{\scalebox{1.5}{$\circ$}}="a3"
        !{(0,0);a(135)**{}?(1)}*{\scalebox{1.5}{$\circ$}}="a4"
        !{(0,0);a(180)**{}?(1)}*{\scalebox{1.5}{$\circ$}}="a5"
        !{(0,0);a(225)**{}?(1)}*{\scalebox{1.5}{$\circ$}}="a6"
        !{(0,0);a(270)**{}?(1)}*{\scalebox{1.5}{$\circ$}}="a7"
        !{(0,0);a(315)**{}?(1)}*{\scalebox{1.5}{$\circ$}}="a8"
        %
        !{(0,0);a(0)**{}?(1.75)}*{\scalebox{1.5}{$\circ$}}="h1"
        !{(0,0);a(45)**{}?(1.75)}*{\scalebox{1.5}{$\circ$}}="h2"
        !{(0,0);a(90)**{}?(1.75)}*{\scalebox{1.5}{$\circ$}}="h3"
        !{(0,0);a(135)**{}?(1.75)}*{\scalebox{1.5}{$\circ$}}="h4"
        !{(0,0);a(180)**{}?(1.75)}*{\scalebox{1.5}{$\circ$}}="h5"
        !{(0,0);a(225)**{}?(1.75)}*{\scalebox{1.5}{$\circ$}}="h6"
        !{(0,0);a(270)**{}?(1.75)}*{\scalebox{1.5}{$\circ$}}="h7"
        !{(0,0);a(315)**{}?(1.75)}*{\scalebox{1.5}{$\circ$}}="h8"
        %
        "a1"-"a2" "a2"-"a3" "a3"-"a4" "a4"-"a5"
        "a5"-"a6" "a6"-"a7" "a7"-"a8" "a8"-"a1"
        "a1"-"h1" "a2"-"h2" "a3"-"h3" "a4"-"h4"
        "a5"-"h5" "a6"-"h6" "a7"-"h7" "a8"-"h8"
      }
    \]
    \vspace{-10 pt}\caption*{\(C_1(C_8\odot K_1)\)}
  \end{minipage}
    \caption{Reduction of $C_3(C_8 \odot K_3)$}
    \label{C reduct-img}
\end{figure}
\begin{theorem} \label{coronathrm}
For the graph \(C_n \odot K_p\), the \(k\)-conversion number is given by
\[
C_k(C_n \odot K_p) =
\begin{cases}
(k-1)n + 1 & \text{if } k \le p + 1, \\
pn + \left\lceil \frac{n}{2} \right\rceil & \text{if } k = p + 2, \\
k\text{-inconvertible} & \text{if } k \ge p + 3.
\end{cases}
\]
\end{theorem}

\begin{proof}
We consider three cases based on the value of \(k\).

\par\medskip
\noindent\textbf{Case 1.} \(k \le p+1\).

By Lemma~\ref{C reduct}, we have:
\[
C_k(C_n \odot K_p) = (k-1)n + C_1(C_n \odot K_{p - (k-1)}).
\]
Let \(p^* = p - (k - 1)\), so \(p^* \ge 0\). Then:
\[
C_k(C_n \odot K_p) = (k - 1)n + C_1(C_n \odot K_{p^*}).
\]
Since the graph \(C_n \odot K_{p^*}\) has no isolated vertices, its 1-conversion number is trivially 1. Hence:
\[
C_k(C_n \odot K_p) = (k - 1)n + 1, \quad \text{for } k \le p+1.
\]

\par\medskip
\noindent\textbf{Case 2.} \(k = p+2\).

Again, by Lemma~\ref{C reduct},
\[
C_k(C_n \odot K_p) = pn + C_2(C_n \odot K_0).
\]
Note that \(K_0\) is the empty graph, so \(C_n \odot K_0 \cong C_n\). Therefore, by Proposition~\ref{prop-cycle},
\[
C_2(C_n \odot K_0) = C_2(C_n) = \left\lceil \frac{n}{2} \right\rceil.
\]
Thus,
\[
C_k(C_n \odot K_p) = pn + \left\lceil \frac{n}{2} \right\rceil, \quad \text{for } k = p+2.
\]

\par\medskip
\noindent\textbf{Case 3.} \(k \ge p+3\).

Applying Lemma~\ref{C reduct} once more,
\[
C_k(C_n \odot K_p) = pn + C_{k^*}(C_n \odot K_0) = pn + C_{k^*}(C_n), \;\;\;\text{where \(k^* = k - p\).}
\]

Since each vertex in \(C_n\) has degree 2, and by Remark~\ref{deg}, \(C_{k^*}(C_n)\) is \(k^*\)-inconvertible for all \(k^* \ge 3\), it follows that:
\[
C_k(C_n \odot K_p) \text{ is } k\text{-inconvertible, for all } k \ge p+3.
\]
\end{proof}
\noindent

\section{Threshold Dynamics on Double Corona Product Graphs}\label{DC_chap}

\subsection{Structural Observations and Boundary Cases}

\begin{lemma}\label{DC.L1}
    Given any $k$-conversion set $S_0$ of $C_n \Dodot K_p$, there exists a related $k$-conversion set $S'_0$ such that $|S'_0|\leq|S_0|$ with each block containing at most one seed vertex.
\end{lemma}
\begin{proof} 
We begin by making several observations. 
\begin{itemize}
\item Any block with two vertices, either seed, converted, or a mixture of both complete the next iteration. As such, any irreversible 2-conversion set with any blocks containing more than 2 seeds is suboptimal and need not be considered in determining the minimal size.
\item \(B\) converts if and only if there exists an iteration such that the two adjacent blocks are complete: \(CBC\). 
\item \(MOMMMMMM\ldots M\) is an irreversible 2-conversion set without any blank blocks and with as many seed vertices as blocks.
\end{itemize}

Let \(S_0\) be an irreversible 2-conversion set with at least one blank block. We can consider \(S_0\) as a collection of strings of non-blank blocks book-ended with blank blocks (possibly the same blank block on both ends). If we can show that each string can be restructured so that the blocks are still converted with at most the same number of seed vertices and having at most one seed vertex in each block, we will have shown the result. We proceed with three cases. Let \(N\) be the number of non-blank blocks between two consecutive blank blocks.

\par\medskip
\noindent\textbf{Case 1.} \(N=1\).

The string must be \(BTB\). If \(BTB\) converts, then both \(OMB\) and \(IMB\) convert.

\par\medskip
\noindent\textbf{Case 2.} \(N=2\).

The string must be \(BTMB\) or \(BTTB\). If either string converts, then both \(OMMB\) and \(IMMB\) convert.

\par\medskip
\noindent\textbf{Case 3.} \(N>2\).

There are at least \(N\) seed vertices. \(BMOM\underbrace{MM\ldots M}_{N-3}B\) will convert with \(N\) seed vertices. 
\end{proof}
\begin{theorem}\label{DC.BC1}
$C_2(C_n \Dodot K_p) = \lfloor \frac{3n+3}{4}\rfloor =n-\lfloor\frac{n}{4}\rfloor, \text{ for }\, p \ge 2$.
\end{theorem}
\begin{proof}
By Lemma~\ref{DC.L1}, we need only consider sets with at most one seed vertex in each block. Suppose that \(S_0\) is such an irreversible 2-conversion set and has more than $\lfloor\frac{n}{4}\rfloor$ blank blocks. By the pidgeon-hole principle, there exists at least one string of four blocks with 2 blank blocks, but this necessitates two seed vertices in a single block which is a contradiction. As such, any irreversible 2-conversion set with at most one seed vertex in each block cannot have more than $\lfloor\frac{n}{4}\rfloor$ many blank blocks meaning that there must be at least $n-\lfloor\frac{n}{4}\rfloor$ many seed vertices.

To establish the result, we provide explicit constructions of such sets of minimal size.
\begin{align*}
n \, \equiv \, 0 \, (\rm mod \,  4)&~~~~\underbrace{MOMB}_1\dots\underbrace{MOMB}_{\frac{n}{4}} \\
n \, \equiv \, 1 \, (\rm mod \,  4)&~~~~~\underbrace{MOMB}_1\dots\underbrace{MOMB}_{\lfloor\frac{n}{4}\rfloor}M \\
n \, \equiv \, 2 \, (\rm mod \,  4)&~~~~~\underbrace{MOMB}_1\dots\underbrace{MOMB}_{\lfloor\frac{n}{4}\rfloor}MO \\
n \, \equiv \, 3 \, (\rm mod \,  4)&~~~~~\underbrace{MOMB}_1\dots\underbrace{MOMB}_{\lfloor\frac{n}{4}\rfloor}MOM
\end{align*}

\end{proof}
\begin{theorem}\label{DC.BC2}
$C_2(C_n \Dodot K_1) = n $.
\end{theorem}

\begin{proof}

By Lemma~\ref{DC.L1}, it is sufficient to consider 2-conversion sets in which each block of $C_n \Dodot K_1$ contains at most one seed vertex. Suppose $S_0$ is an irreversible 2-conversion set with fewer than $n$ seed vertices. 
Since there are $n$ blocks and each block contains at most one seed, at least one block must be blank. Consider such a blank block. Since each adjacent block contains at most one seed, the blank block cannot simultaneously have both its inner and outer vertex be colored from the colored neighbors. Therefore, the blank block cannot become colored which is a contradiction, giving the lower bound:
\[C_2(C_n \Dodot K_1) \ge n\]
A construction using exactly one seed per block of the form $\{I,O,I,O, \cdots, I,O\} \,$ or \\ 
$\{I,O,I,O, \, \cdots, I,O, M \}$ if $n$ is even or odd respectively,  successfully colors the graph, hence 
\[
C_2(C_n \Dodot K_1) \le n.
\]
Therefore, \[C_2(C_n \Dodot K_1) = n\]
\end{proof}

\begin{theorem}\label{DC.BC3}
$ C_2(C_n \Dodot K_0) = 2 \lceil \frac{n}{2}\rceil \,$
\end{theorem}
\begin{proof}
    As a direct result of Remark~\ref{special_bodot}, 
    \[C_2(C_n \Dodot K_0) = C_k\left(\bigcup_{\ell=1}^2 C_n^{[\ell]}\right)= C_k\left(C_n \cup C_n^{\prime}\right)= 2 \cdot C_2(C_n) = 2\left\lceil \frac{n}{2}\right\rceil\]
\end{proof}

\subsection{Reduction Lemma and Threshold Characterization}

\begin{lemma}
    \label{DC reduct}
     $C_k(C_n \Dodot K_p) = r\cdot n + C_{k-r}(C_n \Dodot K_{p-r})$ where $r=\min\{(k-2), p\}$
\end{lemma}

\begin{proof}
    From Remark \ref{deg}, we know that $(k-2)n$ $M$ nodes are necessary in $S_0$ for the graph to be completely colored. As such, every vertex, regardless of whether it is in a block or in a cycle, is connected to at least $k-2$ colored vertices. As such, we may consider the associated graph without these $(k-2)n$ colored vertices and with the lowered threshold number of 2.
\[C_k(C_n \protect\Dodot K_p) = (k-2)n + C_2(C_n \protect\Dodot K_{p-k+2}), \quad \text{ for } k \le p+2\]

On the other hand, if \( k > p+2 \), all \( p \) vertices in each block are removed. The problem then reduces to completing the graph \(C_n \Dodot K_0\) under a threshold of \(k-p \). Hence,
\[C_k(C_n \protect\Dodot K_p) = pn + C_{k-p}(C_n \protect\Dodot K_{p-k+2}), \quad \text{ for } k > p+2\]

It follows that:
\[
C_k(C_n \Dodot K_p) = 
\begin{cases}
(k-2)n + C_{2}(C_n \Dodot K_{p - (k-2)}) & \text{if } k \le p+2,
\\
pn + C_{k - p}(C_n \Dodot K_0) & \text{if } k > p+2.
\end{cases}
\]
\end{proof}

\begin{figure}[H]
    \centering
  \begin{minipage}{0.245\textwidth}
    \centering \[ \xygraph{
!{<0cm,0cm>;<0cm,0.85 cm>:<-0.85 cm,0cm>::}
!{(0,0);a(0)**{}?(1)}*{\scalebox{1.5}{$\circ$}}="a1"
!{(0,0);a(72)**{}?(1)}*{\scalebox{1.5}{$\circ$}}="a2" 
!{(0,0);a(144)**{}?(1)}*{\scalebox{1.5}{$\circ$}}="a3" 
!{(0,0);a(216)**{}?(1)}*{\scalebox{1.5}{$\circ$}}="a4" 
!{(0,0);a(288)**{}?(1)}*{\scalebox{1.5}{$\circ$}}="a5"
\\
!{(0,0);a(0)**{}?(2.5)}*{\scalebox{1.5}{$\circ$}}="c1"
!{(0,0);a(72)**{}?(2.5)}*{\scalebox{1.5}{$\circ$}}="c2" 
!{(0,0);a(144)**{}?(2.5)}*{\scalebox{1.5}{$\circ$}}="c3" 
!{(0,0);a(216)**{}?(2.5)}*{\scalebox{1.5}{$\circ$}}="c4" 
!{(0,0);a(288)**{}?(2.5)}*{\scalebox{1.5}{$\circ$}}="c5"
\\ 
!{(0,0);a(-10)**{}?(1.5)}*{\scalebox{1.5}{$\bullet$}}="b11"
!{(0,0);a(10)**{}?(1.5)}*{\scalebox{1.5}{$\bullet$}}="b12"
!{(0,0);a(0)**{}?(1.9)}*{\scalebox{1.5}{$\circ$}}="b13"
\\ 
!{(0,0);a(62)**{}?(1.5)}*{\scalebox{1.5}{$\bullet$}}="b21"
!{(0,0);a(82)**{}?(1.5)}*{\scalebox{1.5}{$\bullet$}}="b22"
!{(0,0);a(72)**{}?(1.9)}*{\scalebox{1.5}{$\circ$}}="b23"
\\ 
!{(0,0);a(134)**{}?(1.5)}*{\scalebox{1.5}{$\bullet$}}="b31"
!{(0,0);a(154)**{}?(1.5)}*{\scalebox{1.5}{$\bullet$}}="b32"
!{(0,0);a(144)**{}?(1.9)}*{\scalebox{1.5}{$\circ$}}="b33"
\\ 
!{(0,0);a(206)**{}?(1.5)}*{\scalebox{1.5}{$\bullet$}}="b41"
!{(0,0);a(226)**{}?(1.5)}*{\scalebox{1.5}{$\bullet$}}="b42"
!{(0,0);a(216)**{}?(1.9)}*{\scalebox{1.5}{$\circ$}}="b43"
\\ 
!{(0,0);a(278)**{}?(1.5)}*{\scalebox{1.5}{$\bullet$}}="b51"
!{(0,0);a(298)**{}?(1.5)}*{\scalebox{1.5}{$\bullet$}}="b52"
!{(0,0);a(288)**{}?(1.9)}*{\scalebox{1.5}{$\circ$}}="b53"
\\ 
"a1"-"a2"   "a2"-"a3"   "a3"-"a4"   "a4"-"a5"   "a5"-"a1"
\\ 
"c1"-"c2"   "c2"-"c3"   "c3"-"c4"   "c4"-"c5"   "c5"-"c1"
\\ 
"b11"-"b12" "b12"-"b13" "b13"-"b11" 
"b21"-"b22" "b22"-"b23" "b23"-"b21" 
"b31"-"b32" "b32"-"b33" "b33"-"b31" 
"b41"-"b42" "b42"-"b43" "b43"-"b41" 
"b51"-"b52" "b52"-"b53" "b53"-"b51" 
\\ 
"a1"-"b11"  "a1"-"b12"  "a1"-"b13"
"a2"-"b21"  "a2"-"b22"  "a2"-"b23"
"a3"-"b31"  "a3"-"b32"  "a3"-"b33"
"a4"-"b41"  "a4"-"b42"  "a4"-"b43"
"a5"-"b51"  "a5"-"b52"  "a4"-"b43"
\\ 
"c1"-"b11"  "c1"-"b12"  "c1"-"b13" 
"c2"-"b21"  "c2"-"b22"  "c2"-"b23" 
"c3"-"b31"  "c3"-"b32"  "c3"-"b33" 
"c4"-"b41"  "c4"-"b42"  "c4"-"b43" 
"c5"-"b51"  "c5"-"b52"  "c5"-"b53" 
}
\]
\vspace{-10 pt}\caption*{\hspace{28 pt}\(C_4(C_5 \Dodot K_3)\)}
  \end{minipage}
  \hspace{40 pt}
  \begin{minipage}{0.325\textwidth}
    \centering \[ \xygraph{
!{<0cm,0cm>;<0cm,0.85 cm>:<-0.85 cm,0cm>::}
!{(0,0);a(0)**{}?(1)}*{\scalebox{1.5}{$\circ$}}="a1"
!{(0,0);a(72)**{}?(1)}*{\scalebox{1.5}{$\circ$}}="a2" 
!{(0,0);a(144)**{}?(1)}*{\scalebox{1.5}{$\circ$}}="a3" 
!{(0,0);a(216)**{}?(1)}*{\scalebox{1.5}{$\circ$}}="a4" 
!{(0,0);a(288)**{}?(1)}*{\scalebox{1.5}{$\circ$}}="a5"
\\
!{(0,0);a(0)**{}?(2.5)}*{\scalebox{1.5}{$\circ$}}="c1"
!{(0,0);a(72)**{}?(2.5)}*{\scalebox{1.5}{$\circ$}}="c2" 
!{(0,0);a(144)**{}?(2.5)}*{\scalebox{1.5}{$\circ$}}="c3" 
!{(0,0);a(216)**{}?(2.5)}*{\scalebox{1.5}{$\circ$}}="c4" 
!{(0,0);a(288)**{}?(2.5)}*{\scalebox{1.5}{$\circ$}}="c5"
\\ 
!{(0,0);a(-15)**{}?(1.4)}*{\textcolor{light_gray}{\scalebox{1.00}{$\bullet$}}}="b11"
!{(0,0);a(15)**{}?(1.4)}*{\textcolor{light_gray}{\scalebox{1.00}{$\bullet$}}}="b12"
!{(0,0);a(0)**{}?(1.75)}*{\scalebox{1.5}{$\circ$}}="b13"
\\ 
!{(0,0);a(57)**{}?(1.4)}*{\textcolor{light_gray}{\scalebox{1.00}{$\bullet$}}}="b21"
!{(0,0);a(87)**{}?(1.4)}*{\textcolor{light_gray}{\scalebox{1.00}{$\bullet$}}}="b22"
!{(0,0);a(72)**{}?(1.75)}*{\scalebox{1.5}{$\circ$}}="b23"
\\ 
!{(0,0);a(129)**{}?(1.4)}*{\textcolor{light_gray}{\scalebox{1.00}{$\bullet$}}}="b31"
!{(0,0);a(159)**{}?(1.4)}*{\textcolor{light_gray}{\scalebox{1.00}{$\bullet$}}}="b32"
!{(0,0);a(144)**{}?(1.75)}*{\scalebox{1.5}{$\circ$}}="b33"
\\ 
!{(0,0);a(201)**{}?(1.4)}*{\textcolor{light_gray}{\scalebox{1.00}{$\bullet$}}}="b41"
!{(0,0);a(231)**{}?(1.4)}*{\textcolor{light_gray}{\scalebox{1.00}{$\bullet$}}}="b42"
!{(0,0);a(216)**{}?(1.75)}*{\scalebox{1.5}{$\circ$}}="b43"
\\ 
!{(0,0);a(273)**{}?(1.4)}*{\textcolor{light_gray}{\scalebox{1.00}{$\bullet$}}}="b51"
!{(0,0);a(303)**{}?(1.4)}*{\textcolor{light_gray}{\scalebox{1.00}{$\bullet$}}}="b52"
!{(0,0);a(288)**{}?(1.75)}*{\scalebox{1.5}{$\circ$}}="b53"
\\ 
"a1"-"a2"   "a2"-"a3"   "a3"-"a4"   "a4"-"a5"   "a5"-"a1"
\\ 
"c1"-"c2"   "c2"-"c3"   "c3"-"c4"   "c4"-"c5"   "c5"-"c1"
\\ 
"b11"-@{.}"b12" "b12"-@{.}"b13" "b13"-@{.}"b11" 
"b21"-@{.}"b22" "b22"-@{.}"b23" "b23"-@{.}"b21" 
"b31"-@{.}"b32" "b32"-@{.}"b33" "b33"-@{.}"b31" 
"b41"-@{.}"b42" "b42"-@{.}"b43" "b43"-@{.}"b41" 
"b51"-@{.}"b52" "b52"-@{.}"b53" "b53"-@{.}"b51" 
\\ 
"a1"-@{.}"b11"  "a1"-@{.}"b12"  "a1"-"b13"
"a2"-@{.}"b21"  "a2"-@{.}"b22"  "a2"-"b23"
"a3"-@{.}"b31"  "a3"-@{.}"b32"  "a3"-"b33"
"a4"-@{.}"b41"  "a4"-@{.}"b42"  "a4"-"b43"
"a5"-@{.}"b51"  "a5"-@{.}"b52"  "a5"-"b53"
\\ 
"c1"-@{.}"b11"  "c1"-@{.}"b12"  "c1"-"b13"
"c2"-@{.}"b21"  "c2"-@{.}"b22"  "c2"-"b23"
"c3"-@{.}"b31"  "c3"-@{.}"b32"  "c3"-"b33"
"c4"-@{.}"b41"  "c4"-@{.}"b42"  "c4"-"b43"
"c5"-@{.}"b51"  "c5"-@{.}"b52"  "c5"-"b53"
}
\]
\vspace{-10 pt}\caption*{\hspace{0 pt}\(C_{4-2}(C_5 \Dodot K_{3-2})\)}
  \end{minipage}
  \hfill
  \begin{minipage}{0.245\textwidth}
    \centering \[ \xygraph{
!{<0cm,0cm>;<0cm,0.85 cm>:<-0.85 cm,0cm>::}
!{(0,0);a(0)**{}?(1)}*{\scalebox{1.75}{$\circ$}}="a1"
!{(0,0);a(72)**{}?(1)}*{\scalebox{1.75}{$\circ$}}="a2" 
!{(0,0);a(144)**{}?(1)}*{\scalebox{1.75}{$\circ$}}="a3" 
!{(0,0);a(216)**{}?(1)}*{\scalebox{1.75}{$\circ$}}="a4" 
!{(0,0);a(288)**{}?(1)}*{\scalebox{1.75}{$\circ$}}="a5"
\\
!{(0,0);a(0)**{}?(2.5)}*{\scalebox{1.75}{$\circ$}}="c1"
!{(0,0);a(72)**{}?(2.5)}*{\scalebox{1.75}{$\circ$}}="c2" 
!{(0,0);a(144)**{}?(2.5)}*{\scalebox{1.75}{$\circ$}}="c3" 
!{(0,0);a(216)**{}?(2.5)}*{\scalebox{1.75}{$\circ$}}="c4" 
!{(0,0);a(288)**{}?(2.5)}*{\scalebox{1.75}{$\circ$}}="c5"
\\ 
!{(0,0);a(0)**{}?(1.75)}*{\scalebox{1.75}{$\circ$}}="b12"
\\ 
!{(0,0);a(72)**{}?(1.75)}*{\scalebox{1.75}{$\circ$}}="b22"
\\ 
!{(0,0);a(144)**{}?(1.75)}*{\scalebox{1.75}{$\circ$}}="b32"
\\ 
!{(0,0);a(216)**{}?(1.75)}*{\scalebox{1.75}{$\circ$}}="b42"
\\ 
!{(0,0);a(288)**{}?(1.75)}*{\scalebox{1.75}{$\circ$}}="b52"
\\ 
"a1"-"a2"   "a2"-"a3"   "a3"-"a4"   "a4"-"a5"   "a5"-"a1"
\\ 
"c1"-"c2"   "c2"-"c3"   "c3"-"c4"   "c4"-"c5"   "c5"-"c1"
\\ 
\\ 
"a1"-"b12"  
"a2"-"b22"  
"a3"-"b32"  
"a4"-"b42"  
"a5"-"b52" 
\\ 
"c1"-"b12"  
"c2"-"b22"  
"c3"-"b32"  
"c4"-"b42"  
"c5"-"b52" 
}
\]
\vspace{-10 pt}\caption*{\hspace{30 pt}\(C_2(C_5 \Dodot K_1)\)}
\end{minipage}
    \caption{Reduction of $C_4(C_5 \Dodot K_3)$}
    \label{DC reduct-img}   
\end{figure}

\begin{theorem} \label{dcthrm}
For the graph \(C_n \Dodot K_p\), the \(k\)-conversion number is given by
\[
C_k(C_n \Dodot K_p) =
\begin{cases}
        (k-2)n+ \lfloor\frac{3n+3}{4}\rfloor & \text{ if } k \le p
        \\[-0.5 em] 
         pn & \text{ if } k = p+1
        \\[-0.5 em]
        pn + 2\lceil\frac{n}{2}\rceil  & \text{ if } k = p+2
        \\[-0.5 em]
        k\text{-inconvertible} & \text{ if } k \ge p+3
    \end{cases}
\]
\end{theorem}

\begin{proof}
We consider four cases based on the value of \(k\).

\par\medskip
\noindent\textbf{Case 1.} \(k \le p\).

As a direct consequence of Lemma~\ref{DC reduct}, we have:
\[
C_k(C_n \Dodot K_p) = (k-2)n + C_2(C_n \Dodot K_{p - (k-2)}).
\]
Let \(p^* = p - (k - 2)\), so \(p^* \ge 2\). Then:
\[
C_k(C_n \Dodot K_p) = (k - 2)n + C_2(C_n \Dodot K_{p^*}).
\]
By Theorem~\ref{DC.BC1}, 
\[
C_2(C_n \Dodot K_{p^*}) = \left\lfloor\frac{3n+3}{4}\right\rfloor, \quad \text{for } p^* \ge 2
\]
Thus,
\[
C_k(C_n \Dodot K_p) = (k - 2)n + \left\lfloor\frac{3n+3}{4}\right\rfloor, \quad \text{for } k \le p.
\]

\par\medskip
\noindent\textbf{Case 2.} \(k = p+1\).

By Lemma~\ref{DC reduct},
\[
C_k(C_n \Dodot K_p) = (k-2)n + C_2(C_n \Dodot K_1).
\]
By Theorem~\ref{DC.BC2},
\[
C_2(C_n \Dodot K_1) = n.
\]
Therefore,
\[
C_k(C_n \Dodot K_p) = (k-2)n + n = pn, \quad \text{for } k=p+1.
\]

\par\medskip
\noindent\textbf{Case 3.} \(k = p+2\).

Again, by Lemma~\ref{DC reduct},
\[
C_k(C_n \Dodot K_p) = pn + C_2(C_n \Dodot K_0).
\]
From Theorem~\ref{DC.BC3}, we obtain:
\[
C_2(C_n \Dodot K_0) = 2\cdot C_2(C_n) = 2\left\lfloor \frac{n}{2} \right\rfloor.
\]
Thus,
\[
C_k(C_n \Dodot K_p) = pn + 2\left\lfloor \frac{n}{2} \right\rfloor, \quad \text{for } k = p+2.
\]

\par\medskip
\noindent\textbf{Case 4.} \(k \ge p+3\).

Applying Lemma~\ref{DC reduct} once more,
\[
C_k(C_n \Dodot K_p) = pn + C_{k^*}(C_n \Dodot K_0) = pn + C_{k^*}(C_n),
\]
where \(k^* \ge k-p\).

Since each vertex in \(C_n\) has degree 2, and by Remark~\ref{deg}, \(C_{k^*}(C_n)\) is \(k^*\)-inconvertible for all \(k^* \ge 3\), it follows that:
\[
C_k(C_n \Dodot K_p) \text{ is } k\text{-inconvertible, for all } k \ge p+3.
\]
\end{proof}
\noindent

\subsection{Comparison of Corona and Double Corona Threshold Dynamics}\label{comparison}

Although $C_n \odot K_p$ and $C_n \Dodot K_{p-1}$ have the same order and size, their irreversible
$k$-threshold dynamics differ due to fundamental structural differences. In the corona product $C_n \odot K_p$, each attached copy of $K_p$ is connected to the base cycle by a single attachment vertex, so
conversion within a block must propagate through that single bridge.

In contrast, each block of the double corona product $C_n \Dodot K_{p-1}$ contains two attachment vertices (inner and outer), providing two independent bridges between the cycle structure and the complete
graph. This allows a copy of $K_{p-1}$ to become saturated from two sides, significantly altering
the conditions under which vertices convert and enabling configurations that do not occur in
the standard corona product.

These structural differences explain why the threshold behavior of
$C_n \Dodot K_{p-1}$ cannot be obtained by a direct parameter shift from
$C_n \odot K_p$, despite their identical order and size.

\section{Threshold Dynamics on Base-$b$ Corona Product Graphs}
\subsection{Observations of Special Cases}
\begin{remark}\label{remark: disjoint}
By Remark~\ref{special_bodot}, 
\[
  G \bodot{b} K_0 \;\cong\; \bigcup_{\ell=1}^b G^{[\ell]},
\]
Then 
\[
  C_k(G \bodot{b} K_0) = C_k \left( \: \bigcup_{\ell=1}^b G^{[\ell]} \right) = b \cdot C_k (G)
\]
\end{remark}

\begin{remark}\label{rem: p=0}
    For $p=0$, by Proposition~\ref{prop-cycle} and Remark~\ref{remark: disjoint},
    \[C_k(C_n \bodot{b} K_p)= C_k(C_n \bodot{b} K_0) = b \cdot C_k(C_n)=
    \begin{cases}
    b & \text{if } k = 1, \\
    b\lceil \frac{n}{2}\rceil & \text{if } k = 2, \\
    k\text{-inconvertible} & \text{if } k \ge 3.
    \end{cases}
    \]
\end{remark}

\begin{remark}\label{rem: k=1}
When $k=1$ and $p \ge 1$, the graph $C_n \bodot{b} K_p$ has no isolated vertices, hence by Proposition~\ref{k=1 spread},
\[
C_k(C_n \bodot{b} K_p) = 1.
\]
\end{remark}

\subsection{Reduction Lemma and Threshold Characterization of Degenerate Cases}

\begin{proposition}\label{prop:block}
If $p \ge 1$ and $k \ge 2$, every $k$-conversion seed set $S_0$ of $C_n \bodot{b} K_p$ must satisfy
\[\max\{0,\,k-b\}\;\le\;|S_0\cap K_p^{[i]}|\;\le\;p.\]

\end{proposition}

\begin{proof}
We consider two cases.

\par\medskip
\noindent\textbf{Case 1.} \(2\le k\le b\).

Since \(\max\{0,k-b\}=0\), the inequality is trivially satisfied.

\par\medskip
\noindent\textbf{Case 2.} \(b<k\le p+b\).

Before any activation occurs inside \(K_p^{[i]}\), the only colored vertices available within it are those initially contained in \(S_0\cap K_p^{[i]}\).
Therefore, for any vertex in \(K_p^{[i]}\) to be converted at any stage, it is necessary that
\[
|S_0\cap K_p^{[i]}|\;\ge\;k-b = \max\{0, k-b\}.
\]
Since \(K_p^{[i]}\) has exactly \(p\) vertices, we also have \(|S_0\cap K_p^{[i]}|\le p\).
Combining both gives
\[
\max\{0,k-b\}\;\le\;|S_0\cap K_p^{[i]}|\;\le\;p.
\]
At the boundary \(k=p+b\), \(\max\{0, k-b\} = p\), therefore
\[
|S_0\cap K_p^{[i]}|= p.
\]
\end{proof}

\begin{remark}
When $k > p + b$, the quantity $k-b$ exceeds the number of vertices in $K_p^{[i]}$. 
Hence the process can only succeed trivially by coloring all vertices of $K_p^{[i]}$, 
that is \(|S_0 \cap K_p^{[i]}| =  p.\)
\end{remark}

\begin{lemma} \label{lem:generic-reduction}
If $p \ge 1$ and $k\ge 2$,
\[C_k(C_n \bodot{b} K_p) = r \cdot n + C_{k-r} (C_n \bodot{b} K_{p-r}), \quad \text{where } r = \min\{p, \max\{0, k-b\}\}\]

Moreover,
\[
C_k(C_n \bodot{b} K_p) = 
\begin{cases}
 C_k(C_n \bodot{b} K_{p})& \text{if } 2 \le k \le b, \\
 (k-b)n + C_b(C_n \bodot{b} K_{p-k+b})& \text{if } b < k \le p+b, \\
 pn + C_{k-p}(C_n \bodot{b} K_0)& \text{if } k > p+b.
\end{cases}
\]
\end{lemma}

\begin{proof}
Let $r$ be the number of vertices in each $K_p^{[i]}$ that are required in $S_0$. By Proposition~\ref{prop:block}, each $K_p^{[i]}$ must contribute at least $\max\{0,k-b\}$ and at most $p$ seeds to $S_0$, hence $r = \min\{p, \max\{0, k-b\}\}$.
We can then consider the residual graph obtained by removing exactly $r$ colored vertices from each $K_p^{[i]}$ for all $i$, a total of $r \cdot n$ vertices. The resulting graph $C_n \bodot{b} K_{p-r}$ then requires only a $(k-r)$-threshold for completion, as each remaining vertex is adjacent to at least $r$ colored vertices.

Hence for $r = \min\{p, \max\{0, k-b\}\}$, 
\[
C_k(C_n \bodot{b} K_p) = r \cdot n + C_{k-r} (C_n \bodot{b} K_{p-r})
\]
In particular, we have
\[
r=
\begin{cases}
 0& \text{if } 2 \le k \le b, \\
 k-b& \text{if } b < k \le p+b, \\
 p& \text{if } k > p+b.
\end{cases}
\]
which, upon substitution, yields the stated piecewise result.
\end{proof}

\begin{figure}[H]
    \centering
    \begin{minipage}{0.4\textwidth}
\includegraphics[height=7cm,width=\linewidth]{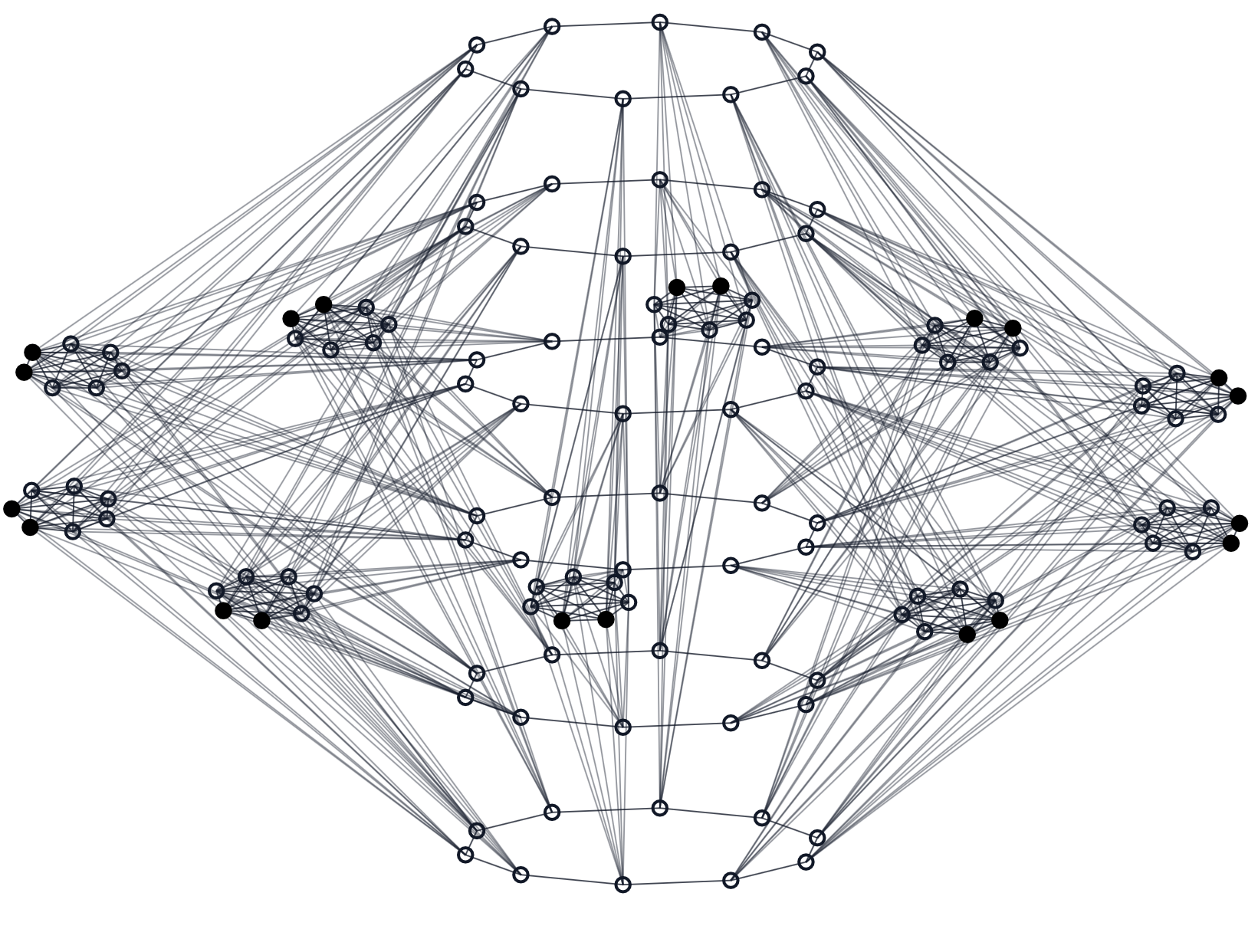}
\vspace{-25 pt}
    \caption*{\(C_{8}(C_{10} \bodot{6} K_7)\)}
    \end{minipage}
    \hfill
    \begin{minipage}{0.4925\textwidth}
\includegraphics[height=7cm,width=\linewidth]{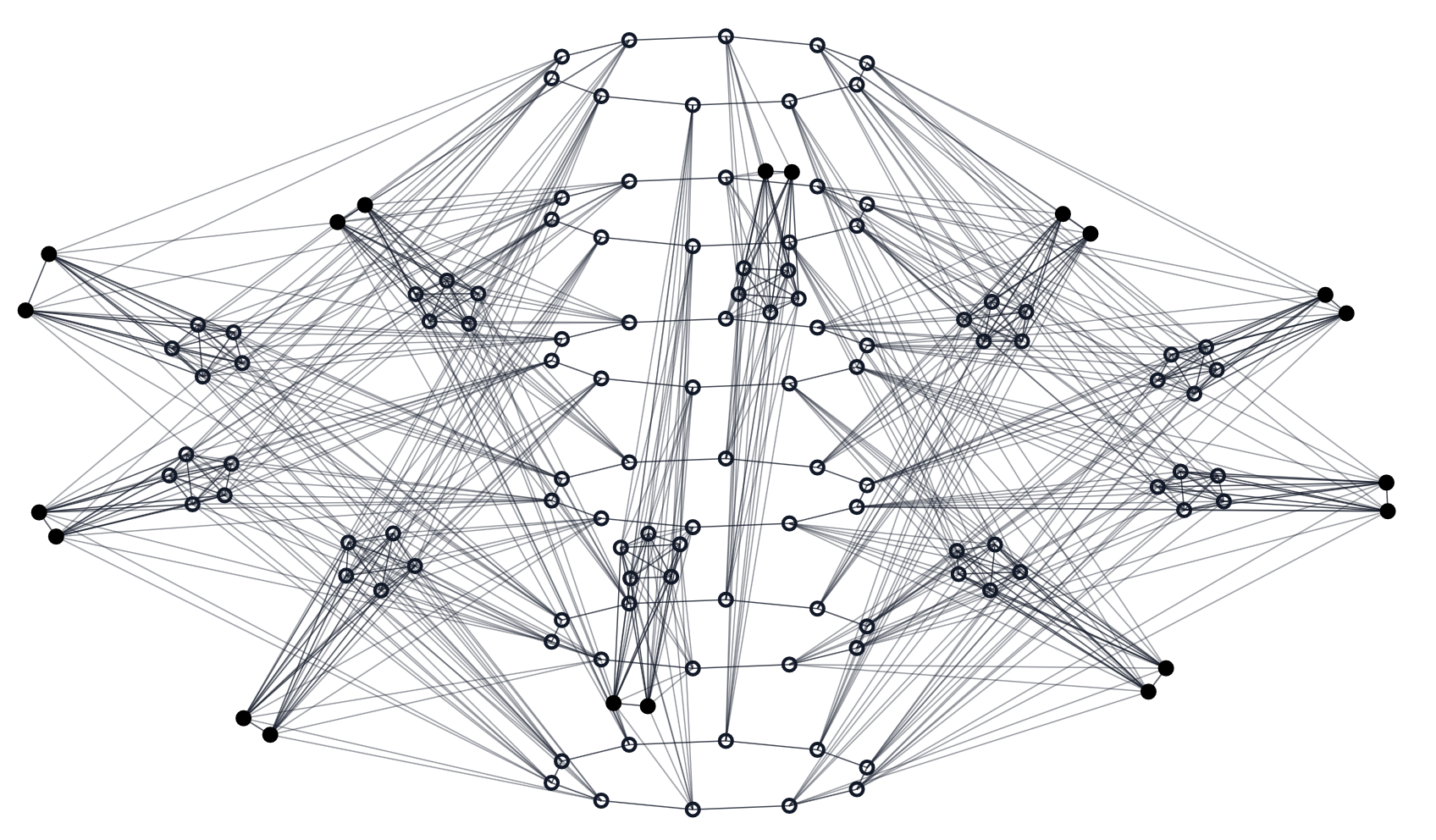}
\vspace{-20 pt}
\caption*{\(C_{8-2}(C_{10} \bodot{6} K_{7-2})\)}
    \end{minipage}
    \hfill
    \begin{minipage}{0.45\textwidth}
\includegraphics[height=7cm,width=\linewidth]{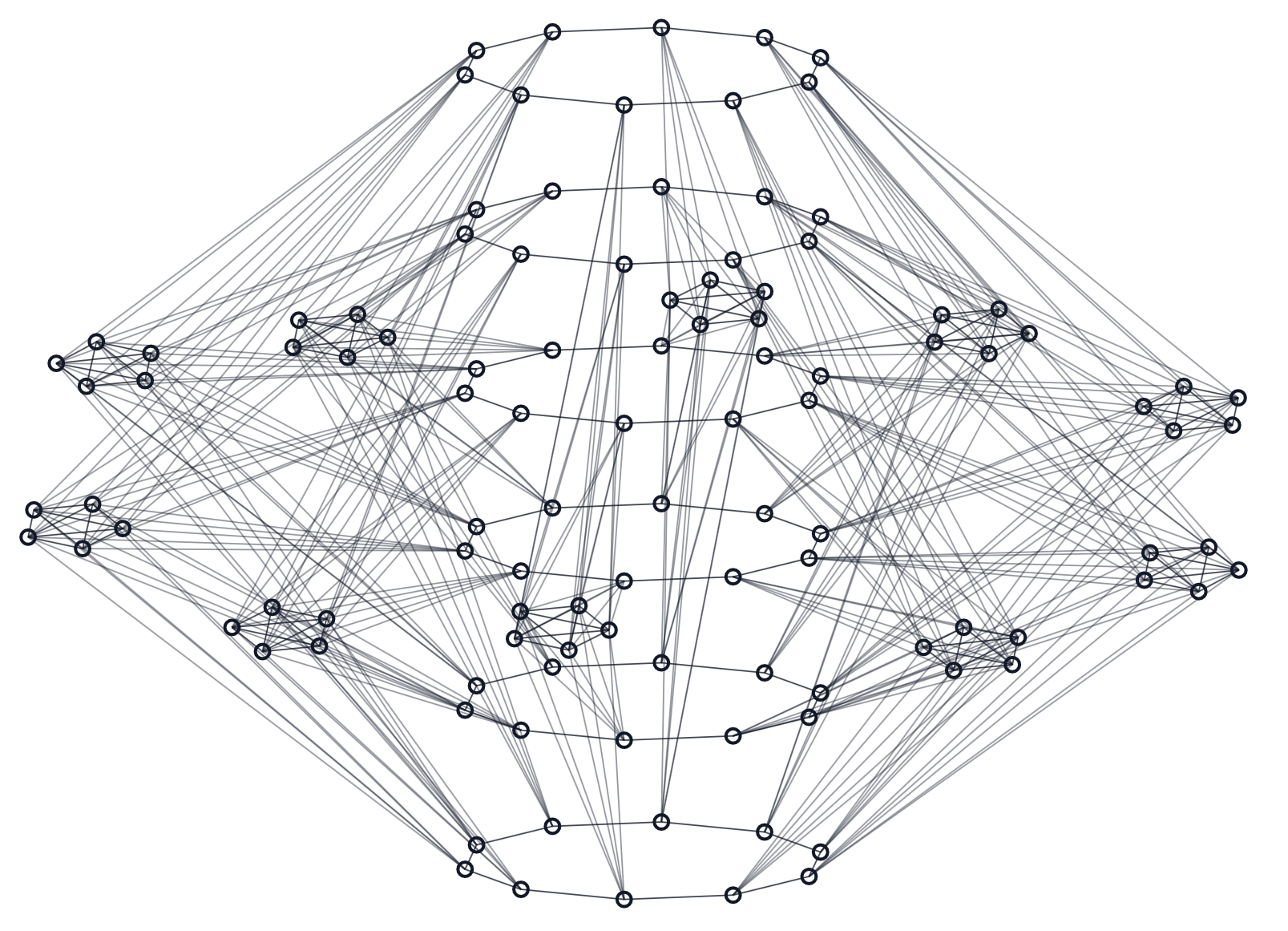}
\vspace{-20 pt}
\caption*{\(C_{6}(C_{10} \bodot{6} K_5)\)}
    \end{minipage}
    \caption{Reduction of \(C_{8}(C_{10} \bodot{6} K_7)\)}
\end{figure}

\noindent The following corollary summarizes a general lower and upper bound implied by Lemma~\ref{lem:generic-reduction}.

\begin{corollary}\label{k-bound}
Let $r=\min\{p,\max\{0,k-b\}\}$. For all $k\ge2$,
\[
 rn + \gamma \le C_k(C_n\bodot{b} K_p) \le n(p+b),\qquad
\gamma=\begin{cases}
1,& r\le p-1,\\[2pt]
b,& r=p.
\end{cases}
\]
The lower bound is tight when $b=1$ with $1<k\le p+1$ and also $k=p+1$ for any $b\ge 1$.
\end{corollary}

\begin{proof}
From Lemma~\ref{lem:generic-reduction},
\(C_k(C_n\bodot{b} K_p)=rn+C_{k-r}(C_n\bodot{b}K_{p-r}).\)
Since $k\ge2$ and $r=\min\{p,\max\{0,k-b\}\}$, it follows that $k-r\ge1$ in all cases; hence the residual graph always requires at least one seed.
\begin{itemize}
    \item If $p-r\ge1$, the residual graph is is connected. By Proposition~\ref{k=1 spread},\newline $C_{k-r}(C_n \bodot{b} K_{p-r})\ge 1$.  
    \item If instead $p-r=0$, then by Remark~\ref{rem: p=0}, \(C_{k-r}(C_n\bodot{b}K_{0})=b\cdot C_{k-r}(C_n)\ge b.\)
\end{itemize}
This establishes the lower bound. Equality occurs when the residual term attains its minimum value, that is,
\[
C_{k-r}(C_n\bodot{b}K_{p-r})=
\begin{cases}
1,& r=p-1,\\[2pt]
b,& r=p,
\end{cases}
\]
yielding
\[
C_k(C_n\bodot{b} K_p)=
\begin{cases}
(k-b)n+1,& \text{if } b=1 \text{ and } 1<k < p+1,\\[2pt]
pn+b,& \text{if } k=p+1.
\end{cases}
\]
\\
The upper bound \(C_k(C_n \bodot{b} K_p)\le n(p+b) = |V(C_n \bodot{b} K_p)|\) is trivial.

\end{proof}

\begin{proposition}\label{prop:high-threshold}
If $k > p+b$, then
\[
C_k(C_n \bodot{b} K_p) =
\begin{cases}
pn + b\lceil \frac{n}{2}\rceil, & \text{if } k = p+2 \text{ and } b=1, \\[0.25em]
k\text{-inconvertible}, & \text{if } k \ge p+3 \text{ and } b \ge 2.
\end{cases}
\]
\end{proposition}
\begin{proof}
As a direct result of Lemma~\ref{lem:generic-reduction},
\[
C_k(C_n \bodot{b} K_p) = pn + C_{k-p}(C_n \bodot{b} K_0)
\quad \text{when } k > p+b.
\]

By Remark~\ref{remark: disjoint},
\[
C_{k-p}(C_n \bodot{b} K_0) = b \cdot C_{k-p}(C_n).
\]

Let $k^* := k-p$. From Proposition~\ref{prop-cycle},
\[
C_{k^*}(C_n) =
\begin{cases}
1, & \text{if } k^* = 1, \\
\lceil \frac{n}{2}\rceil, & \text{if } k^* = 2, \\
k^*\text{-inconvertible}, & \text{if } k^* \ge 3.
\end{cases}
\]

Since $k>p+b$, we have $k^*>b\ge1$. In particular:
\begin{itemize}
    \item The case $k^*=1$ (equivalently $k=p+1$) cannot occur.
    \item The case $k^*=2$ occurs only when $b=1$.
    \item If $b\ge2$, then $k^*>b$ implies $k^*\ge3$, and each residual cycle layer is $k^*$-inconvertible.
\end{itemize}

Therefore,
\begin{itemize}
    \item If $k=p+2$ and $b=1$, then
    \[
    C_k(C_n \bodot{b} K_p)=pn+b\lceil\tfrac{n}{2}\rceil \implies C_k(C_n \bodot{1} K_p) = pn +\lceil\tfrac{n}{2}\rceil. 
    \]
    as seen in Chapter \ref{C Chap}.
    \item If $k\ge p+3$ and $b\ge2$, the residual graph is $k^*$-inconvertible, and hence $C_n \bodot{b} K_p$ is $k$-inconvertible.
\end{itemize}

Thus, except for the boundary case $b=1$, the high-threshold condition $k>p+b$ forces $C_n \bodot{b} K_p$ to be $k$-inconvertible.
\end{proof}

\section{Probabilistic Analysis}\label{probability}
\subsection{Probabilistic $k$-Conversion and Resistance Factors}
Classical results on irreversible $k$-threshold processes focus on the conversion number
$C_k(G)$, the minimum size of a seed set that could result in full conversion of the graph.
While this deterministic viewpoint identifies extremal thresholds, it does not address how
likely a randomly chosen seed set of that size is to succeed.

To address this limitation, we introduce a probabilistic refinement of the threshold framework.
In addition to determining the size $C_k(G)$, we ask how frequently
such sets occur among all subsets of that size.

For a graph $G$ and integer $k \ge 1$, let
\[
\mathcal{M}_k(G) := \{ S_0 \subseteq V(G) : |S_0| = C_k(G),\; S_0 \text{ is an I$k$CS} \}.
\]

We define the \emph{irreversible $k$-conversion probability} of $G$ by
\[
\mathbf{P}_{k}\,[\;\! G \;\!] := \frac{|\mathcal{M}_k(G)|}{\dbinom{|V(G)|}{C_k(G)}} \, ,
\]  

Intuitively, $\mathbf{P}_{k}[G]$ measures the likelihood that a uniformly random choice of $C_k(G)$ vertices
forms a valid irreversible $k$-conversion seed set. This probability depends not only on the
threshold value $C_k(G)$, but also on the structural distribution of successful seed sets in $G$.

These quantities also enable comparative analysis of threshold dynamics across different graph families. 

\subsection{Probabilistic Threshold Dynamics on Classical Graph Families}

In this chapter, we study probabilistic irreversible $k$-threshold dynamics on several graph families. We begin with classical graphs, where minimum irreversible $k$-conversion sets and their success probabilities can be described directly from local degree constraints and symmetry. The structure of corona product graphs introduces new combinatorial behavior that requires a more refined probabilistic analysis.

\begin{remark}\label{prob_k1}
If $k=1$ and $G$ is connected, then $C_1(G)=1$ and every minimum initial set consists of a single vertex whose color propagates to the entire graph. Hence $\mathbf P_1[G]=1$.
\end{remark}

\begin{remark}\label{prob_k-inconvertible}
If $G$ is $k$-inconvertible, then $C_k(G)=|V(G)|$ and the only minimum initial
set is $S_0=V(G)$, yielding $\mathbf P_k[G]=1$.
\end{remark}

\subsubsection{Complete Graphs $K_n$}
In a complete graph, every vertex is adjacent to all others. Consequently, if at least $k$ vertices are initially colored, every uncolored vertex has at least $k$ colored neighbors and converts immediately. Hence,
\[
C_k(K_n)=k \quad \text{for } 1\le k\le n,
\]
and for $k>n$ the graph is $k$-inconvertible. In either case, every minimum irreversible $k$-conversion set succeeds. Therefore,
\[
\mathbf P_k[K_n]=1 \quad \text{for all } k\ge 1.
\]

\subsubsection{Paths $P_n$}

For the path graph $P_n$, it is known \cite{dreyer2009} that
\(
C_2(P_n)=\left\lceil \frac{n+1}{2}\right\rceil.
\)
This arises from the boundary structure of $P_n$, the endpoints have degree one and must belong to any irreversible $2$-conversion set. Once the endpoints are fixed, the remaining vertices must alternate to ensure that each interior vertex eventually has two colored neighbors.

When $n$ is even, there are exactly two minimum $2$-conversion sets, corresponding to the two possible alternating patterns. When $n$ is odd, the alternating pattern is unique. Consequently,
\[
\mathbf P_2[P_n]=\frac{2-[n\!\!\!\!\pmod 2]}{\binom{n}{\lceil \frac{n+1}{2}\rceil}}.
\]

\begin{table}[H]
    \centering
    \renewcommand{\arraystretch}{1.3}
    \setlength{\tabcolsep}{10pt}
    \begin{tabular}{ccccc}
    \toprule
    $n$ & 5 & 8 & 14 & 37 \\
    \midrule
    $\mathbf{P}_2[P_n]$ & $0.1000$ & $0.0357$ & $6.6600\times 10^{-4}$ & $5.6585\times 10^{-11}$ \\
    \bottomrule
    \end{tabular}
\caption{Sample values of $\mathbf{P}_2[P_n]$ illustrating rapid decay with $n$}
    \label{tab:path-prob-scale}
\end{table}

\subsubsection{Cycles $C_n$}\label{cycle-prob}

Unlike paths, cycles have no boundary vertices. For the cycle graph $C_n$, irreversible $2$-conversion sets must alternate around the cycle in order for every vertex to eventually receive two colored neighbors. The structure of such alternating configurations depends on the parity of $n$.

If $n$ is even, there are exactly two minimum $2$-conversion sets, corresponding to selecting either all odd-indexed or all even-indexed vertices. If $n$ is odd, any minimum $2$-conversion set must contain exactly one pair of consecutive vertices, with the remaining vertices alternating around the cycle. Since the cycle is rotationally symmetric, there are $n$ such
configurations.

\begin{figure}[H] 
\begin{minipage}{0.20\textwidth}
\end{minipage}
\hfill
  \begin{minipage}{0.15\textwidth}
    \centering \[ \xygraph{
!{<0cm,0cm>;<0cm,0.55 cm>:<0.55 cm,0cm>::}
!{(0,0);a(0)**{}?(1)}*{\scalebox{1.25}{$\bullet$}}="a1" 
        !{(0,0);a(0)**{}?(1.75)}*{v_1}
!{(0,0);a(60)**{}?(1)}*{\scalebox{1.25}{$\circ$}}="a2" 
        !{(0,0);a(60)**{}?(1.75)}*{v_2}
!{(0,0);a(120)**{}?(1)}*{\scalebox{1.25}{$\bullet$}}="a3" 
        !{(0,0);a(120)**{}?(1.75)}*{v_3}
!{(0,0);a(180)**{}?(1)}*{\scalebox{1.25}{$\circ$}}="a4"
        !{(0,0);a(180)**{}?(1.75)}*{v_4}
!{(0,0);a(240)**{}?(1)}*{\scalebox{1.25}{$\bullet$}}="a5"
        !{(0,0);a(240)**{}?(1.75)}*{v_5}
!{(0,0);a(300)**{}?(1)}*{\scalebox{1.25}{$\circ$}}="a6"
        !{(0,0);a(300)**{}?(1.75)}*{v_6}
\\ 
"a1"-"a2"   "a2"-"a3"   "a3"-"a4"   "a4"-"a5"   "a5"-"a6"   "a6"-"a1"
\\
    !{(0,0);a(180)**{}?(2.5)}*{\text{\scalebox{0.9}{Without loss of}}}
    !{(0,0);a(180)**{}?(3.25)}*{\text{\scalebox{0.9}{generality, choose}}}
    !{(0,0);a(180)**{}?(4)}*{\text{\scalebox{0.9}{odd-indexed vertices}}}
}
\]
  \end{minipage}
  \hfill
  \begin{minipage}{0.025\textwidth}
  \centering \[ \xygraph{
    !{<0cm,0cm>;<0cm,0.55cm>:<0.55cm,0cm>::}
    !{(0,0);a(0)**{}?(1)}*{\scalebox{1.25}{ }}="a1"
    !{(0,0);a(90)**{}?(1)}*{\scalebox{1.25}{ }}="a2"
    !{(0,0);a(180)**{}?(1)}*{\scalebox{1.25}{ }}="a3"
    \\
            !{(0,0);a(180)**{}?(0)}*{\scalebox{0.75}{$\times$} \, 2}
            \\
            \\
            !{(0,0);a(180)**{}?(1.8)}*{\scalebox{0.65}{\text{Unique}}}
            !{(0,0);a(180)**{}?(2.3)}*{\scalebox{0.65}{\text{Rotational}}}
            !{(0,0);a(180)**{}?(2.9)}*{\scalebox{0.65}{\text{ Symmetry}}}
    !{(0,0);a(270)**{}?(1)}*{\scalebox{1.25}{ }}="a4"
    \\
    "a1" :@/^0.2cm/ "a2"  "a3" :@/^0.2cm/ "a4"
  }
  \]
\end{minipage}
  \hfill
  \begin{minipage}{0.15\textwidth}
    \centering \[ \xygraph{
!{<0cm,0cm>;<0cm,0.55 cm>:<0.55 cm,0cm>::}
!{(0,0);a(0)**{}?(1)}*{\scalebox{1.25}{$\circ$}}="a1" 
        !{(0,0);a(0)**{}?(1.75)}*{v_1}
!{(0,0);a(60)**{}?(1)}*{\scalebox{1.25}{$\bullet$}}="a2" 
        !{(0,0);a(60)**{}?(1.75)}*{v_2}
!{(0,0);a(120)**{}?(1)}*{\scalebox{1.25}{$\circ$}}="a3" 
        !{(0,0);a(120)**{}?(1.75)}*{v_3}
!{(0,0);a(180)**{}?(1)}*{\scalebox{1.25}{$\bullet$}}="a4"
        !{(0,0);a(180)**{}?(1.75)}*{v_4}
!{(0,0);a(240)**{}?(1)}*{\scalebox{1.25}{$\circ$}}="a5"
        !{(0,0);a(240)**{}?(1.75)}*{v_5}
!{(0,0);a(300)**{}?(1)}*{\scalebox{1.25}{$\bullet$}}="a6"
        !{(0,0);a(300)**{}?(1.75)}*{v_6}
\\ 
"a1"-"a2"   "a2"-"a3"   "a3"-"a4"   "a4"-"a5"   "a5"-"a6"   "a6"-"a1"
\\
        !{(0,0);a(180)**{}?(2.5)}*{\text{\scalebox{0.9}{even-indexed vertices}}}
}
\]
  \end{minipage}
  \hfill
  \begin{minipage}{0.20\textwidth}
\end{minipage}
\hfill
  \newline
  \begin{minipage}{0.20\textwidth}
    \centering \[ \xygraph{
!{<0cm,0cm>;<0cm,0.45 cm>:<0.45 cm,0cm>::}
!{(0,0);a(0)**{}?(1)}*{\scalebox{1.25}{$\bullet$}}="a1" 
        !{(0,0);a(0)**{}?(1.75)}*{v_1}
!{(0,0);a(72)**{}?(1)}*{\scalebox{1.25}{$\circ$}}="a2" 
        !{(0,0);a(72)**{}?(1.75)}*{v_2}
!{(0,0);a(144)**{}?(1)}*{\scalebox{1.25}{$\bullet$}}="a3" 
        !{(0,0);a(144)**{}?(1.75)}*{v_3}
!{(0,0);a(216)**{}?(1)}*{\scalebox{1.25}{$\bullet$}}="a4"
        !{(0,0);a(216)**{}?(1.75)}*{v_4}
!{(0,0);a(288)**{}?(1)}*{\scalebox{1.25}{$\circ$}}="a5"
        !{(0,0);a(288)**{}?(1.75)}*{v_5}
\\
        !{(0,0);a(180)**{}?(2.5)}*{\text{\scalebox{0.9}{Without loss of}}}
        !{(0,0);a(180)**{}?(3.25)}*{\text{\scalebox{0.9}{generality, choose}}}
        !{(0,0);a(180)**{}?(4)}*{\text{\scalebox{0.9}{\(S_0= \{v_1, v_3, v_4\}\)}}}
\\ 
"a1"-"a2"   "a2"-"a3"   "a3"-"a4"   "a4"-"a5"   "a5"-"a1"
}
\]
  \end{minipage}
  \hspace{5 pt}
  \begin{minipage}{0.125\textwidth}
  \centering \[ \hspace{-10 pt} \xygraph{
    !{<0cm,0cm>;<0cm,0.55cm>:<0.55cm,0cm>::}
    !{(0,0);a(0)**{}?(1)}*{\scalebox{1.25}{ }}="a1"
    !{(0,0);a(90)**{}?(1)}*{\scalebox{1.25}{ }}="a2"
    !{(0,0);a(180)**{}?(1)}*{\scalebox{1.25}{ }}="a3"
    !{(0,0);a(270)**{}?(1)}*{\scalebox{1.25}{ }}="a4"
    \\      
            !{(0,0);a(180)**{}?(0)}*{\scalebox{0.75}{$\times$} \, n}
            \\
            !{(0,0);a(180)**{}?(1.8)}*{\scalebox{0.65}{\text{Unique}}}
            !{(0,0);a(180)**{}?(2.3)}*{\scalebox{0.65}{\text{Rotational}}}
            !{(0,0);a(180)**{}?(2.9)}*{\scalebox{0.65}{\text{ Symmetry}}}
    \\
    "a1" :@/^0.2cm/ "a2"  "a3" :@/^0.2cm/ "a4"
  }
  \]
\end{minipage}
  \hfill    
  \begin{minipage}{0.15\textwidth}
    \centering \[\xygraph{
!{<0cm,0cm>;<0cm,0.45 cm>:<0.45 cm,0cm>::}
!{(0,0);a(0)**{}?(1)}*{\scalebox{1.25}{$\circ$}}="a1" 
        !{(0,0);a(0)**{}?(1.75)}*{v_1}
!{(0,0);a(72)**{}?(1)}*{\scalebox{1.25}{$\bullet$}}="a2" 
        !{(0,0);a(72)**{}?(1.75)}*{v_2}
!{(0,0);a(144)**{}?(1)}*{\scalebox{1.25}{$\circ$}}="a3" 
        !{(0,0);a(144)**{}?(1.75)}*{v_3}
!{(0,0);a(216)**{}?(1)}*{\scalebox{1.25}{$\bullet$}}="a4"
        !{(0,0);a(216)**{}?(1.75)}*{v_4}
!{(0,0);a(288)**{}?(1)}*{\scalebox{1.25}{$\bullet$}}="a5"
        !{(0,0);a(288)**{}?(1.75)}*{v_5}
\\ 
"a1"-"a2"   "a2"-"a3"   "a3"-"a4"   "a4"-"a5"   "a5"-"a1"
}
\]


  \end{minipage}
  \hfill
  \begin{minipage}{0.15\textwidth}
    \centering \[ \xygraph{
!{<0cm,0cm>;<0cm,0.45 cm>:<0.45 cm,0cm>::}
!{(0,0);a(0)**{}?(1)}*{\scalebox{1.25}{$\bullet$}}="a1" 
        !{(0,0);a(0)**{}?(1.75)}*{v_1}
!{(0,0);a(72)**{}?(1)}*{\scalebox{1.25}{$\circ$}}="a2" 
        !{(0,0);a(72)**{}?(1.75)}*{v_2}
!{(0,0);a(144)**{}?(1)}*{\scalebox{1.25}{$\bullet$}}="a3" 
        !{(0,0);a(144)**{}?(1.75)}*{v_3}
!{(0,0);a(216)**{}?(1)}*{\scalebox{1.25}{$\circ$}}="a4"
        !{(0,0);a(216)**{}?(1.75)}*{v_4}
!{(0,0);a(288)**{}?(1)}*{\scalebox{1.25}{$\bullet$}}="a5"
        !{(0,0);a(288)**{}?(1.75)}*{v_5}
\\ 
"a1"-"a2"   "a2"-"a3"   "a3"-"a4"   "a4"-"a5"   "a5"-"a1"
}
\]
  \end{minipage}
  \hfill
    \begin{minipage}{0.15\textwidth}
    \centering \[ \xygraph{
!{<0cm,0cm>;<0cm,0.45 cm>:<0.45 cm,0cm>::}
!{(0,0);a(0)**{}?(1)}*{\scalebox{1.25}{$\bullet$}}="a1" 
        !{(0,0);a(0)**{}?(1.75)}*{v_1}
!{(0,0);a(72)**{}?(1)}*{\scalebox{1.25}{$\bullet$}}="a2" 
        !{(0,0);a(72)**{}?(1.75)}*{v_2}
!{(0,0);a(144)**{}?(1)}*{\scalebox{1.25}{$\circ$}}="a3" 
        !{(0,0);a(144)**{}?(1.75)}*{v_3}
!{(0,0);a(216)**{}?(1)}*{\scalebox{1.25}{$\bullet$}}="a4"
        !{(0,0);a(216)**{}?(1.75)}*{v_4}
!{(0,0);a(288)**{}?(1)}*{\scalebox{1.25}{$\circ$}}="a5"
        !{(0,0);a(288)**{}?(1.75)}*{v_5}
\\ 
"a1"-"a2"   "a2"-"a3"   "a3"-"a4"   "a4"-"a5"   "a5"-"a1"
}
\]
  \end{minipage}
  \hfill
  \begin{minipage}{0.15\textwidth}
    \centering \[ \xygraph{
!{<0cm,0cm>;<0cm,0.45 cm>:<0.45 cm,0cm>::}
!{(0,0);a(0)**{}?(1)}*{\scalebox{1.25}{$\circ$}}="a1" 
        !{(0,0);a(0)**{}?(1.75)}*{v_1}
!{(0,0);a(72)**{}?(1)}*{\scalebox{1.25}{$\bullet$}}="a2" 
        !{(0,0);a(72)**{}?(1.75)}*{v_2}
!{(0,0);a(144)**{}?(1)}*{\scalebox{1.25}{$\bullet$}}="a3" 
        !{(0,0);a(144)**{}?(1.75)}*{v_3}
!{(0,0);a(216)**{}?(1)}*{\scalebox{1.25}{$\circ$}}="a4"
        !{(0,0);a(216)**{}?(1.75)}*{v_4}
!{(0,0);a(288)**{}?(1)}*{\scalebox{1.25}{$\bullet$}}="a5"
        !{(0,0);a(288)**{}?(1.75)}*{v_5}
\\ 
"a1"-"a2"   "a2"-"a3"   "a3"-"a4"   "a4"-"a5"   "a5"-"a1"
}
\]
  \end{minipage}
  \hfill
    \caption{The number of unique minimum irreversible $2$-conversion sets for a cycle is 2 (or $n$) when $n$ is even (or odd).}
    \label{prob-img}
\end{figure}

Thus, the number of distinct alternating $2$-conversion sets on $C_n$ is
\[
2 + (n-2)[n \!\!\!\!\pmod 2],
\]
hence
\[
\mathbf P_2[C_n]=
\frac{2 + (n-2)[n \!\!\!\!\pmod 2]}{\dbinom{n}{\lceil \frac{n}{2}\rceil}}.
\]

\begin{table}[H]
    \centering
    \renewcommand{\arraystretch}{1.3}
    \setlength{\tabcolsep}{10pt}
    \begin{tabular}{ccccc}
    \toprule
    $n$ & 5 & 8 & 14 & 37 \\
    \midrule
    $\mathbf{P}_2[C_n]$ & $0.5000$ & $0.0286$ & $5.8275\times 10^{-4}$ & $2.0936\times 10^{-9}$ \\
    \bottomrule
    \end{tabular}
\caption{Sample values of $\mathbf{P}_2[C_n]$ illustrating rapid decay with $n$}
    \label{tab:path-prob-scale}
\end{table}

\subsubsection{Complete Bipartite Graphs $K_{m,n}$}
We now analyze probabilistic threshold dynamics on the complete bipartite graph $K_{m,n}$, whose vertex set is partitioned into two independent sets of sizes $m$ and $n$.

It is known \cite{dreyer2009} that the irreversible $k$-conversion number of the complete bipartite graph $K_{m,n}$ satisfies

\[
C_k(K_{m,n}) =
\begin{cases}
k & \text{if }~\min\{m,n\} \ge k 
\\
\max\{m,n\} & \text{if }~ \max\{m,n\} \ge k > \min\{m,n\}
\\
k\text{-inconvertible} & \text{if }~ k > \max\{m,n\}.
\end{cases}
\]
\\[-2.0 em]

Since vertices have neighbors only in the opposite partite set, a vertex can convert under the irreversible $k$-threshold rule only if at least $k$ vertices in the opposite part are initially colored. This structural constraint leads to several cases depending on the value of $k$.

\par\medskip
\noindent\textbf{Case 1. } {$2 \le k \le \min\{m,n\}$} \\
A vertex in either partite set has at least $k$ neighbors available in the opposite part. For conversion to occur, all initial seeds must therefore be contained entirely within a single partite set. Hence, distributing seeds across both sides fails to provide any vertex with enough colored neighbors to trigger conversion. Therefore, a successful initial set consists of choosing $k$ vertices from one side or the other. The number of successful configurations is
\begin{align}
    \binom{m}{k} + \binom{n}{k}.
\end{align}
Since the total number of minimum initial conversion sets of size $k$ is $\binom{m+n}{k}$, the successful probability is
\begin{align}
    \mathbf P_k[K_{m,n}] = \dfrac{\binom{m}{k} + \binom{n}{k}}{\binom{m+n}{k}}
\end{align}

\par\medskip
\noindent\textbf{Case 2. } {$\min\{m,n\} < k \le \max\{m,n\}$} \\
In this case, vertices in the larger partite set (of size $\max\{m,n\}$) have degree $\min\{m,n\}$, which is strictly less than $k$ and can never convert unless they are initially colored. Consequently, any successful initial set must include all vertices of the larger partite set ($\max\{m,n\}$). 

Thus, there is exactly one successful initial configuration, consisting of the entire larger partite set. Since the minimum irreversible $k$-conversion set has size $\max\{m,n\}$, the probability of success is
\begin{align}
    \mathbf P_k[K_{m,n}] = \dfrac{1}{\binom{m+n}{\max\{m,n\}}}
\end{align}

\begin{table}[H]
    \centering
    \renewcommand{\arraystretch}{1.25}
    \setlength{\tabcolsep}{8pt}
    \begin{tabular}{ccccc}
    \toprule
    $m$ & $n$ & $k$ & Case & $\mathbf{P}_k[K_{m,n}]$ \\
    \midrule
    3 & 5 & 2 & $2 \le k \le \min\{m,n\}$ & $0.2143$  \\
    
    8 & 7 & 4 & $2 \le k \le \min\{m,n\}$ & $0.1026$ \\
    
    10 & 4 & 6 & $\min\{m,n\} < k \le \max\{m,n\}$ & $9.9900\times10^{-4}$ \\
    
    9 & 21 & 13 & $\min\{m,n\} < k \le \max\{m,n\}$ & $6.9895\times10^{-8}$ \\
    \bottomrule
    \end{tabular}
    \caption{Sample probabilities for $K_{m,n}$ graphs, illustrating dependence on $m$, $n$, and $k$.}
    \label{tab:kmn-prob-samples}
\end{table}

\subsection{Probabilistic Threshold Dynamics on $C_n \odot K_p$}
Unlike classical graph families, corona products exhibit layered attachment structures that allow multiple distinct minimum irreversible $k$-conversion configurations. We now analyze how this structure affects probabilistic threshold dynamics by computing exact success probabilities for minimum initial sets for all values of $k$, $n$, and $p$.

\par\medskip
\noindent\textbf{Case 1. }  \(2 \le k \le p + 1\)\\
Let us consider the necessary conditions for a completed graph coloring, by first examining the cliques. Naturally for $k\le p$, \(k-1\) vertices are needed in each of the \(n\) cliques. Then there is one vertex that can be either in one of the cliques or the cycle. The number of ways we can choose $k-1$ vertices in every block and one in the cycle is:

\begin{align}
    n\binom{p}{k-1}^n 
\end{align}
The number of ways we can choose $k-1$ vertices in $n-1$ cliques and $k$ vertices in one of the $n$ cliques is:
\begin{align}
   n\binom{p}{k-1}^{n-1} \binom{p}{k}
\end{align}
Since either one of these options would yield a successful graph, we know that the total number of successful initial coloring is:
\begin{align}
     &n\binom{p}{k-1}^n  + n\binom{p}{k-1}^{n-1} \binom{p}{k} \\
     \notag \\
     &= \quad n\binom{p}{k-1}^{n-1} \binom{p+1}{k}
\end{align}

Thus, the probability of successful colorings is
\begin{align}
   \mathbf{P}_{k}[C_n \odot K_p] = \frac{n\dbinom{p}{k-1}^{n-1} \dbinom{p+1}{k}}{\dbinom{n(p+1)}{(k-1)n+1} } \label{p1}
\end{align}

\par\medskip
\noindent\textbf{Case 2. }  \(k = p + 2\) \\
From Lemma~\ref{C reduct}, we know that we will need to start with all $p$ vertices for each of the $n$ cliques to have a successfully colored graph. The number of ways to do this is:
\begin{align}
    \binom{p}{p}^n = 1
\end{align}

Recall from Lemma~\ref{C reduct}, when $k=p+2$ the conversion dynamics reduce to selecting an alternating $2$-conversion set on the underlying cycle. The number of such configurations is given by the enumeration for $C_n$ derived in Section~\ref{cycle-prob} to be
\begin{align}
2+ (n-2)[n\!\!\!\!\!\pmod 2] 
\end{align}

Thus the probability of successfully selecting initial vertices is
\begin{align}
    \mathbf{P}_{k}[C_n \odot K_p] = \frac{ 2+ (n-2)[n\!\!\!\!\!\pmod2]}{\dbinom{n(p+1)}{pn+\lceil{}{\frac{n}{2}}\rceil{}}}   \label{p2}
\end{align}

\par\medskip
\noindent\textbf{Case 3. }  \(k \ge p + 3\)\\
From Lemma~\ref{C reduct} we found that \(C_n \odot K_p\) is \(k\)-inconvertible, then by Remark~\ref{prob_k-inconvertible}
\begin{align}
\mathbf{P}_{k}[C_n \odot K_p] = 1  \label{p3}
\end{align}

Combining all three above cases, we propose the following theorem. 

\begin{theorem}
For the family of graphs $C_n \odot K_p$, the probability $\mathbf{P}_{k}[C_n \odot K_p]$ that a randomly chosen initial set $S_0$ with $|S_0| = C_k(C_n \odot K_p)$ (the $k$-threshold number) results in a complete $k$-coloring of the graph is given by:
\begin{align*}
        \mathbf{P}_{k}[\,C_n \odot K_p\,] = 
        \begin{cases}
            \dfrac{n\dbinom{p}{k-1}^{n-1}\dbinom{p+1}{k}} {\dbinom{n(p+1)}{(k-1)n+1}}  & \text{if} \quad 2 \le k \le p+1
            \\[2.5 em]
           \dfrac{ 2+ (n-2)[n \!\!\!\!\!\pmod 2]}{\dbinom{n(p+1)}{pn+\lceil{}{\frac{n}{2}}\rceil{}}} & \text{if} \quad k = p+2 
            \\[0.75 em]
            1 & \text{if} \quad k\ge p+3 \\
        \end{cases} \label{ps}
    \end{align*}
\end{theorem}

\begin{table}[H]
    \centering
    \renewcommand{\arraystretch}{1.25}
    \setlength{\tabcolsep}{8pt}
    \begin{tabular}{ccccc}
    \toprule
    $n$ & $p$ & $k$ & Case & $\mathbf{P}_k[C_n\odot K_p]$ \\
    \midrule
    4 & 3 & 2 & $k\le p+1$ & $0.1484$ \\
    
    6 & 11 & 12 & $k\le p+1$ & $4.2883\times10^{-7}$ \\
    
    14 & 7 & 9 & $k=p+2$ & $5.5206\times10^{-11}$ \\
    
    18 & 6 & 8 & $k=p+2$ & $1.2146\times10^{-13}$ \\
    \bottomrule
    \end{tabular}
    \caption{Sample probabilities for $C_n \odot K_p$, illustrating dependence on $n$, $p$, $k$.}
    \label{tab:corona-prob-samples}
\end{table}

\subsection{Critical Cardinality for Guaranteed Success on $C_n \odot K_p$}

Previous results characterize the minimum initial set size
$C_k(C_n \odot K_p)$ for which $\mathbf{P}_k[C_n \odot K_p] > 0$.
Here we determine the smallest cardinality of an initial set that guarantees
successful conversion for \emph{every} possible placement of vertices.

\begin{definition}
The \emph{critical cardinality} $|S_0|_{\mathrm{crit}}$ is the smallest integer such that every initial set $S_0 \subseteq V(C_n \odot K_p)$ with $|S_0| = |S_0|_{\mathrm{crit}}$ satisfies $\mathbf{P}_k[C_n \odot K_p] = 1$ for a fixed $k$.
\end{definition}

Thus, while $\mathbf{P}_k$ captures typical behavior under random placement, $|S_0|_{\mathrm{crit}}$ represents a worst-case threshold for guaranteed success.

For notational convenience, let $v = |V(C_n \odot K_p)| = n(p+1)$ denote the total number of vertices of $C_n \odot K_p$. We consider cases according to the value of $k$.

\par\medskip
\noindent\textbf{Case 1. } {$k \le p + 1$} \\
By Theorem~\ref{coronathrm}, conversion requires at least $k-1$ initially colored vertices in each attached $K_p$ clique. Thus, any failing configuration must contain a clique with at most $k-2$ colored vertices.

To construct the largest such configuration, begin with all $v$ vertices colored and remove as many vertices as possible from a single clique while leaving exactly $k-2$ vertices in that clique. This removes $p-(k-2)$ vertices and yields a failing set of size

\begin{align}
    v - (p-(k-2)) = v - p + k - 2.
\end{align}

Adding one additional vertex eliminates all failing configurations. Hence,
\begin{align}
    |S_0|_{\mathrm{crit}} = v - p + k - 1,
\end{align}
and every initial set of this size succeeds.

\begin{remark}
It follows that $|S_0|_{\mathrm{crit}} \le v-1$ for $k \le p$, while $|S_0|_{\mathrm{crit}} = v$ when $k = p+1$.
\end{remark}

\par\medskip
\noindent\textbf{Case 2. } {$k = p + 2$} \\
Theorem~\ref{coronathrm} implies that each $K_p$ clique must be fully colored initially. Removing even a single vertex produces a clique with only $p-1$ colored vertices, which fails to meet the threshold. Therefore,
\begin{align}
    |S_0|_{\mathrm{crit}} = v.
\end{align}

\par\medskip
\noindent\textbf{Case 3. } {$k > p + 2$} \\
In this case, $C_n \odot K_p$ is $k$-inconvertible unless all vertices are initially colored. Thus,
\begin{align}
    |S_0|_{\mathrm{crit}} = v.
\end{align}

\begin{remark}
The critical cardinality $|S_0|_{\mathrm{crit}}$ marks the point at which $\mathbf{P}_k[C_n \odot K_p]=1$, equivalently when the resistance factor necessarily vanishes. 

If $|S_0|_{\mathrm{crit}} < v$, then the difference $v - |S_0|_{\mathrm{crit}}$ quantifies the structural redundancy of the graph, or the number of vertices that may remain uncolored without risking failure.

If $|S_0|_{\mathrm{crit}} = v$, the graph is \emph{stable}, in the sense that guaranteed conversion is impossible unless all vertices are initially colored.
\end{remark}

\section{Additional Probabilistic Exploration}\label{explore}

Unlike the corona product, the double corona product introduces interactions between two base graph layers that substantially complicate the enumeration of minimum irreversible $k$-conversion sets. While the deterministic threshold behavior of double corona products can be completely characterized, the probabilistic analysis provides no uniform closed-form solution.

Rather than presenting a complete probabilistic theory, this chapter documents exact results of special cases, together with partial progress and structural obstacles that arise in the remaining cases. These results serve both to illustrate the increased combinatorial complexity of double corona products and to motivate future analytical approaches.

\subsection{Exact Probabilistic Results for Special Cases of $C_n~\!\Dodot\;\!K_p$}
For certain $k$-values, the probabilistic behavior of the irreversible $k$-threshold process on $C_n \Dodot K_p$ can be determined exactly.

\par\medskip
\noindent\textbf{Case 1. } {$p>0$, $k=1$}\\
Since $C_n\Dodot K_p$ is connected, the $1$-threshold process converts the graph from any single initial vertex. By Remark~\ref{prob_k1},

\begin{align}
\mathbf{P}_1[C_n \Dodot K_p] = 1.
\end{align}

\par\medskip
\noindent\textbf{Case 2. } {$p=0$, $k=1$}\\
Recall from Remark \ref{special_bodot} that $C_n\Dodot K_0 \cong C_n \cup C_n^{\prime}$. Hence the $1$-threshold process must convert on both of the graphs simultaneously. Since there are $n$ options for each cycle we have, 
\begin{align}
\mathbf{P}_1[C_n \Dodot K_0] = \frac{n^2}{\dbinom{2n}{2}}.
\end{align}

\par\medskip
\noindent\textbf{Case 3. } {$k=p+2$}\\
By Lemma~\ref{DC reduct} (the double corona reduction lemma), once all clique vertices are initially colored, the irreversible $k$-threshold dynamics on the remaining graph decomposes into two independent $2$-threshold processes on the two underlying cycle layers.

Consequently, a successful initial set is obtained by independently selecting a minimum irreversible $2$-conversion set on each of the two disjoint cycle copies. By the enumeration in Section~\ref{cycle-prob}, the number of such configurations on a single cycle $C_n$ is
\begin{align}
    2 + (n-2)[n \!\!\!\!\pmod 2].
\end{align}

Since the two cycle layers evolve independently, the total number of successful initial configurations is therefore
\begin{align}
\Bigl(2 + (n-2)[n \!\!\!\!\pmod 2]\Bigr)^2 .
\end{align}

Thus, the probability of successful conversion in this case is
\begin{align}
\mathbf{P}_{k}[C_n \Dodot K_p]
=\dfrac{\Bigl(2 + (n-2)[n\!\!\!\!\pmod 2]\Bigr)^2}
{\dbinom{n(p+2)}{\,pn + 2\left\lceil \frac{n}{2}\right\rceil\,}} .
\end{align}

\par\medskip
\noindent\textbf{Case 4. } {$k \ge p+3$}\\
By Theorem \ref{dcthrm}, $C_k(C_n \Dodot K_p)$ is $k$-inconvertible. Therefore by Remark \ref{prob_k-inconvertible}, 
\begin{align}
\mathbf{P}_{k}[C_n \Dodot K_p] = 1.
\end{align}

\subsection{Exploratory Probabilistic Analysis of Double Corona Products}
For the remaining threshold values, the probabilistic enumeration of minimum irreversible $k$-conversion sets on $C_n \text{\;\!} \Dodot \text{\:\!} K_p$ becomes substantially more complex. In particular, the case $k = p+1$ exhibits strong interdependence between the two cycle layers and the attached cliques, preventing a clean decomposition into independent subproblems.

In this case, successful configurations depend on how vertices are distributed across both cycle layers and cliques, leading naturally to a partition-based counting framework. While this approach yields partial structural insight, a complete closed-form expression remains open.

Accordingly, we restrict attention to outlining the combinatorial framework and documenting the obstacles encountered, rather than presenting an explicit probability formula. The remaining case $k \le p$ is not pursued here, as it introduces additional layers of interaction beyond those already present in the case $k = p+1$, including overlapping block dependencies across both cycle layers.

\subsubsection{Fixing the Number of Blank Blocks}

Recall from Lemma \ref{DC reduct}, there must be $k-2=p-1$ vertices in each clique, hence a factor of
\begin{align}
    \dbinom{p}{p-1}^n \, =\, p^n
\end{align}

Also from the deterministic analysis, the number of $B$ (blank) blocks is a key factor in whether full conversion can occur. For this reason, we begin the probabilistic analysis by fixing the number of blank blocks, say $x$.

We consider the simplest case, \(x=0\), where no blank blocks are introduced. In this setting, the seed pattern is highly restricted: the inner and outer seed choices must alternate from block to block.

When \(n\) is even, this alternation closes consistently around the cycle, so there are exactly two alternating patterns 
\begin{align}
I,O,I,O,\cdots I,O \qquad\text{or}\qquad O,I,O,I,\cdots O,I     
\end{align}

When \(n\) is odd, a perfect alternation cannot close around the cycle. In that case, one additional $M$ block must bridge the gap between the two cycle layers. Thus, there are $n-1$ blocks with $p-1$ vertices in the cliques and one full block (with exactly $p$ vertices required). By the rotational symmetry of the odd cycle, this $M$ block has $n$ possible placements. Therefore a factor of
\begin{align}
    n\cdot\dbinom{p}{p-1}^{n-1}\cdot\dbinom{p}{p} \,= \,n\,p^{n-1}
\end{align}

Thus,
\begin{align}
\mathbf{P}_{k}[C_n \Dodot K_p] \ge \dfrac{2p^{n-1}\left(p+(n-p)[n\!\!\!\!\pmod2]\right)}{\dbinom{n(p+2)}{pn}}
\end{align}

This gives an explicit lower bound arising from the subclass of successful configurations with no blank blocks. For \(x>0\), the interaction between blank blocks, full blocks, and the two cycle layers introduces additional dependencies that do not admit a comparable closed-form enumeration in the present framework.

Accordingly, the \(x=0\) analysis is recorded here as a first exact contribution in the double corona setting, while the general enumeration problem for \(x>0\) is deferred.

\section{Conclusion}

Irreversible \(k\)-threshold dynamics were studied on corona-type graph products, with emphasis on the corona product, the double corona product, and the base-\(b\) corona product. For \(C_n \odot K_p\) and \(C_n \Dodot K_p\), exact deterministic threshold results were obtained through reduction lemmas that convert the problem to smaller corona-type instances and, ultimately, to classical base graphs. For the base-\(b\) construction, the same reduction perspective yields a unified framework extending the corona and double corona cases.

A probabilistic refinement was also considered. Exact probability formulas were obtained for several classical graph families and for corona product graphs, together with exact special-case results for selected corona-type constructions. These results show that probabilistic behavior can distinguish graph families whose deterministic threshold numbers are closely related, and they illustrate how layered attachment structure influences conversion reliability.

For double corona products, the deterministic theory is complete in the cases treated here, whereas the probabilistic theory is only partially resolved. The analysis of configurations with no blank blocks gives an explicit lower bound in the case \(k=p+1\), but a full enumeration requires additional control over the interaction between blank blocks, full blocks, and the two cycle layers. This identifies a natural direction for subsequent work without affecting the completeness of the deterministic results and the exact probabilistic results established elsewhere in the paper.

\end{document}